\documentclass{article}

\usepackage[hyperfootnotes=false]{hyperref}
\usepackage{notation}
\usepackage{amssymb,bm,amsthm}
\usepackage{color}
\usepackage{todonotes}
\usepackage[symbol]{footmisc}

\DeclareMathOperator{\ex}{\mathbb E}

\DeclareMathOperator{\Vol}{Vol}
\DeclareMathOperator{\tr}{tr}

\DeclareMathOperator{\diag}{diag}
\let\top\intercal

\newtheorem{thm}{Theorem}
\newtheorem{lem}{Lemma}

\newtheorem{cor}{Corollary}

\newtheorem{defi}{Definition}

\newcommand{\KL}{\operatorname{KL}}

\newcommand{\sign}{\operatorname{sign}}
\newcommand{\TV}{\operatorname{TV}}

\newcommand{\eps}{\varepsilon}

\newcommand{\wh}{\widehat}
\newcommand{\wt}{\widetilde}

\newcommand{\ol}{\overline}

\newcommand{\Var}{\operatorname{Var}}
\newcommand{\Cov}{\operatorname{Cov}}

\newcommand{\E}{\mathbb{E}}

\newcommand{\Unif}{\operatorname{Unif}}

\newcommand{\Bias}{\operatorname{Bias}}

\newcommand{\Tr}{\operatorname{tr}}

\newcommand{\ba}{\mathbf{a}}
\newcommand{\bb}{\mathbf{b}}

\newcommand{\be}{\mathbf{e}}

\newcommand{\mB}{\mathcal{B}}

\newcommand{\mF}{\mathcal{F}}
\newcommand{\mG}{\mathcal{G}}
\newcommand{\mM}{\mathcal{M}}
\newcommand{\mN}{\mathcal{N}}

\newcommand{\mV}{\mathcal{V}}

\newcommand{\bA}{\mathbf{A}}
\newcommand{\bB}{\mathbf{B}}
\newcommand{\bC}{\mathbf{C}}
\renewcommand{\ba}{\mathbf{a}}
\newcommand{\bD}{\mathbf{D}}

\newcommand{\bu}{\mathbf{u}}
\newcommand{\bU}{\mathbf{U}}
\newcommand{\bv}{\mathbf{v}}

\newcommand{\bw}{\mathbf{w}}
\newcommand{\bW}{\mathbf{W}}
\newcommand{\bx}{\mathbf{x}}
\newcommand{\bX}{\mathbf{X}}
\newcommand{\by}{\mathbf{y}}

\newcommand{\bxi}{\bm{\xi}}

\newcommand{\btheta}{\bm{\theta}}
\newcommand{\bzeta}{\bm{\zeta}}

\begin{document}

\title{Hebbian learning inspired estimation of the linear regression parameters from queries}

\author{%
    Johannes Schmidt-Hieber\footnotemark[1]%
    \ \ and Wouter M. Koolen\footnotemark[2]%
}
  \footnotetext[1]{University of Twente} 
  \footnotetext[2]{University of Twente and CWI}

\maketitle

\begin{abstract}
Local learning rules in biological neural networks (BNNs) are commonly referred to as Hebbian learning. \cite{2023arXiv230111777S} links a biologically motivated Hebbian learning rule to a specific zeroth-order optimization method. In this work, we study a variation of this Hebbian learning rule to recover the regression vector in the linear regression model. Zeroth-order optimization methods are known to converge with suboptimal rate for large parameter dimension compared to first-order methods like gradient descent, and are therefore thought to be in general inferior. By establishing upper and lower bounds, we show, however, that such methods achieve near-optimal rates if only queries of the linear regression loss are available. Moreover, we prove that this Hebbian learning rule can achieve considerably faster rates than any non-adaptive method that selects the queries independently of the data.
\end{abstract}

%
%
\paragraph{Keywords:} Biological neural networks, derivative-free methods, linear regression model, minimax estimation, sequential estimation, zeroth-order optimization.

\textbf{MSC 2020:} Primary: 62L20; secondary: 62J05

\section{Introduction}

While considerable efforts have been dedicated to develop the theoretical foundations underlying artificial neural networks (ANNs), a theory for learning in biological neural networks (BNNs) remains largely unexplored.

ANNs were designed to mimic BNNs but there are distinct differences in terms of the network structure and the learning. ANNs are deterministic and pass real numbers through the network whereas BNNs are stochastic and biological neurons send so-called spike trains. The information in BNNs is decoded in the spike times, the moments when the neurons fire/discharge. The brain updates the weights locally based on the state of the neurons that it connects. There are various quantitative and qualitative updating rules and they are commonly referred to as Hebbian learning. Specifically, spike-time dependent plasticity rules are Hebbian learning rules based on local spike times. It has been widely acknowledged in previous work that Hebbian learning in BNNs cannot perform gradient descent \cite{Crick1989, Lillicrap2020, FundamentalsCompNeuro}. The rationale behind is that updating one parameter in a gradient descent scheme requires to evaluate a partial derivative that depends on the other parameters. This cannot be implemented by a local learning rule. 

The recent work \cite{2023arXiv230111777S} derives a specific derivative-free optimization method that captures key characteristics of the spike time dependent plasticity rule underlying Hebbian learning in BNNs. For mathematical tractability, we consider here a slight variation of this scheme. 

We study this iterative scheme in a scenario where instead of the full data, one can only query the linear regression model, see Section \ref{sec.query_model} for a definition and discussion. Based on $k$ rounds of querying, the aim of the method is to recover the $d$-dimensional regression vector. Estimation in this setting is non-trivial. In particular, the data are not informative enough to compute gradients and run gradient descent. A first contribution of this article is to derive a bound for the convergence rate of the biologically inspired gradient-free learning rule in the query model. It is moreover shown that up to a log-factor in the number of parameters $d$, the derived convergence rates matches the lower bound for sample sizes $k\gtrsim d^2\log(d)$. The derived minimax lower bound is non-standard. The main obstacle is that sequential estimation procedures with queries depending on previous observations induces dependence among the data. As we do not constraint the queries, the induced dependence on the future data is hard to characterize and to control. Compared to the earlier lower bounds for adaptive sensing  \cite{6289365} and derivative-free stochastic optimization \cite{pmlr-v30-Shamir13, Duchi2015}, the main difficulty here is that the dependence on previous queries and the parameter appears not in the mean but in the variance of the data distribution.

Any method achieving the optimal rate in the linear model with queries needs to carefully exploit the information obtained from previous queries. Indeed, we show that naive methods that specify the queries independently of the data will converge with a slower rate in the dimension parameter $d$, see Corollary \ref{cor.gap} for the precise statement. This should be contrasted with adaptive sensing where it is known that the adaptation to previously seen data can improve the rate by, at most, a log-factor \cite{6289365}. 

The article is structured as follows. Section \ref{sec.main} provides more details on biological learning and the link with gradient-free optimization. After formally stating the query model, the convergence rates are given. A corresponding lower bound is stated in Section \ref{sec.lb_adaptive}. Section \ref{sec.nonadapt} derives matching upper and lower bounds for the convergence rate in the case of non-adaptive queries. This enables us to quantify the gain in the convergence rate by integrating previously observed query values into the query strategy. Related literature is summarized in Section \ref{sec.related_lit}. Proofs are deferred to the Appendix.

{\bf Notation:} For vectors the spectral norm coincides with the Euclidean norm, see Lemma \ref{lem.spect_norm_of_vec_prod}. Therefore, we can denote both norms by $\|\cdot\|.$ Matrix inequalities are taken with respective to the partial ordering of symmetric matrices (Loewner order).

\section{Upper bounds}
\label{sec.main}

\subsection{Hebbian learning inspired update rule}

Working in the supervised learning framework, suppose we want to learn a $d$-dimensional parameter vector $\btheta$ from training data $(\bX_1,Y_1),(\bX_2,Y_2),\ldots,$ 
with $d$-dimensional design vectors/inputs $\bX_1,\bX_2,\ldots$ and real-valued response variables $Y_1,Y_2,\ldots$ For our approach, it is sufficient to assume that the number of parameters and the number of covariates are the same, that is, $\btheta$ and $\bX_1,\bX_2,\ldots$ are all vectors of length $d.$

Whereas gradient descent based methods are ubiquitous in machine learning, BNNs rely on local updating rules and receive feedback through neurotransmitters such as dopamine. Anticipating to solve a task well reduces the amount of released neurotransmitter in the brain. In reward-based synaptic plasticity, the amount of released neurotransmitter is modelled as the difference between a realized loss based on the current task and the anticipated loss that predicts the current loss based on previously seen losses \cite{Fremaux13326}.

Zeroth-order optimization methods are derivative-free methods that use the evaluation of the loss but do not involve the gradient of the loss. Since weight updates in BNNs are based on evaluations of the loss, it seems natural to interpret the learning in BNNs as a specific zeroth-order method. \cite{2023arXiv230111777S} extracts from the spike-time dependent local updating of the weights in a BNN a global zeroth-order rule for learning the parameter vector $\btheta$. If $\btheta_k$ denotes the $k$-th update, $L(\btheta,\bX_k,Y_k)$ is the loss for parameter $\btheta$ on the $k$-th training sample $(\bX_k,Y_k),$ $\ol L_k$ is the anticipated loss in the $k$-th round based on previously seen losses, and $\bU_k$ is a $d$-dimensional uniform random vector $\Unif([A-,A]^d)$, the biologically inspired update formula is
\begin{align}
    \btheta_k 
    &=\btheta_{k-1}+\alpha_k \big( L(\btheta_{k-1}+ \bU_k,\bX_k,Y_k)- \ol L_k \big)\big(e^{-\bU_k}-e^{\bU_k}\big).
    \label{eq.updating}
\end{align}
The expressions $e^{\pm \bU_k}$ have to be understood componentwise, that is, if $\bU_k=(U_{k1},\ldots,U_{kd})^\top,$ then $e^{\bU_k}=(e^{U_{k1}},\ldots, e^{U_{kd}})^\top$ and $e^{-\bU_k}=(e^{-U_{k1}},\ldots, e^{-U_{kd}})^\top.$ The parameter $A$ in the uniform distribution is a constant of the biological network but treated here as a hyperparameter of the optimization method. There is little guidance from the neuroscience on how to choose the anticipated loss $\ol L_k$. In \cite{2023arXiv230111777S} the $\ol L_k$ was taken as the loss from the previous round $L(\btheta_{k-2}+ \bU_{k-1}, \bX_{k-1},Y_{k-1}).$ For mathematical tractability, we assign instead the value $L(\btheta_{k-1}+ \bU_k',\bX_k,Y_k)$ to $\ol L_k,$ for an independently sampled and uniformly distributed random vector $\bU_k'\sim \Unif([-A,A]^d).$ Thus, we will study in this work convergence of the zeroth-order method
\begin{align}
    \btheta_k 
    &=\btheta_{k-1}+\alpha_k \big( L(\btheta_{k-1}+ \bU_k,\bX_k,Y_k)- L(\btheta_{k-1}+ \bU_k',\bX_k,Y_k) \big)\big(e^{-\bU_k}-e^{\bU_k}\big), \quad k=1,\ldots
    \label{eq.updating_special}
\end{align}

Theorem 1 in \cite{2023arXiv230111777S} shows that in expectation, this rule does approximately gradient descent. Working here with a slightly different anticipated loss does not change this result. Indeed, since $\bU_k$ is independent of all other randomness and $\E[e^{-\bU_k}-e^{\bU_k}]=0,$ conditioning on everything except for $\bU_k$ gives
\begin{align*}
    \E\Big[L(\btheta_{k-1}+ \bU_k',\bX_k,Y_k) \big(e^{-\bU_k}-e^{\bU_k}\big)\Big]= 
    \E\Big[L(\btheta_{k-1}+ \bU_k',\bX_k,Y_k) \E\big[e^{-\bU_k}-e^{\bU_k}\big]\Big]
    =0.
\end{align*}
Thus, as in the original version, the anticipated loss vanishes in expectation and Theorem 1 in \cite{2023arXiv230111777S} holds without any changes. While the connection to gradient descent is appealing, the main obstacle for fast convergence of zeroth-order methods is the high variance.

\subsection{The query model}
\label{sec.query_model}

A natural first step of a statistical analysis of biologically inspired learning is to study the properties of the zeroth-order update rule \eqref{eq.updating_special} for standard statistical models. The overall aim is to identify relevant models where this rule achieves the optimal convergence rates and/or outperforms other standard methods. 

A fundamental problem in statistics and machine learning is the case where we want to learn a regression/feature vector $\btheta^\star$ such that for an independent draw $(\bX,Y)$ with the same distribution as the (training) data $(\bX_k,Y_k),$ we have $Y\approx \bX^\top \btheta^\star.$ Considering squared loss 
\begin{align*}
    L(\btheta,\bX,Y)
    = \big(Y - \bX^\top \btheta\big)^2,
\end{align*}
the update formula \eqref{eq.updating_special} can then be rewritten as
\begin{align}
  \btheta_k
  ~=~
  \btheta_{k-1} + \alpha_k \del*{
    \big(Y_k - \bX_k^\top (\btheta_{k-1} + \bU_k)\big)^2
    - \big(Y_k - \bX_k^\top (\btheta_{k-1} + \bU_k')\big)^2
  } \big(e^{-\bU_k}-e^{\bU_k}\big).
  \label{eq.update_with_squared_loss}
\end{align}
Interestingly, this formula does not require full knowledge of the covariate vector $\bX_k$ and can be computed from queries. Here a query is defined as follows. Based on earlier seen queried observations and possibly extra randomness, the statistician chooses in the $k$-th round a query vector $\bv_k.$ Instead of the full $k$-th training sample $(\bX_k,Y_k)$, one can only observe
\begin{align*}
    Z_k =Y_k- \bX_k^\top \bv_k.
\end{align*}
To realize \eqref{eq.update_with_squared_loss}, one needs to query each $(\bX_k,Y_k)$ twice. This means that in the $k$-th round, we choose two query vectors $\bv_k,\bv_k'$ based on previously seen query values and extra randomness. As observations, we get
\begin{align}
    Z_k =Y_k- \bX_k^\top \bv_k, \quad \text{and} \quad Z_k' =Y_k- \bX_k^\top \bv_k'.
    \label{eq.queries_def}
\end{align}
To see that this observational scheme is indeed enough to compute the updates \eqref{eq.update_with_squared_loss}, we can argue by induction, assuming that the initialization $\btheta_0$ has been chosen independently of the data. The induction assumption is then, that $\btheta_k$ only depends on the queries $(Z_1,Z_1',\ldots,Z_k,Z_k')$ and exogenous randomness. Then, $\bv_{k+1}=\btheta_k+\bU_k$ and $\bv_{k+1}'=\btheta_k+\bU_k'$ are admissible query vectors with corresponding queries $Z_{k+1}=Y_{k+1}- \bX_{k+1}^\top(\btheta_k+\bU_{k+1})$ and $Z_{k+1}'=Y_{k+1}- \bX_{k+1}^\top(\btheta_k+\bU_{k+1}').$ The update formula implies that $\btheta_{k+1}$ only depends on the queries $(Z_1,Z_1',\ldots,Z_{k+1},Z_{k+1}')$ and exogenous randomness, completing the induction step.

Queries can be thought of as an attention mechanism that tell us where to look next to extract useful information about the data. A query model seems appropriate if the full input vector cannot be processed because, for instance, the data arrive too quickly. The information from the queries is insufficient to compute gradients. Hence, gradient descent methods cannot be run in this case.

There is a wide variety of previously considered query models in the statistical literature. Examples are querying large graphs to avoid storage problems \cite{JMLR:v22:20-1137} or learning causal relationships from path queries \cite{NEURIPS2018_a0b45d1b}.

To finish the description of the model, we have to specify the distribution of the latent/unobserved variables $(\bX_1,Y_1),(\bX_2,Y_2),\ldots$ We assume that those are i.i.d.\ and generated from the linear regression model with random design, that is, $\bX_1,\bX_2,\ldots $ are drawn i.i.d.\ from some unknown distribution $P_{\bX}$ and the response variables are
\begin{align}
	Y_i= \bX_i^\top \btheta^\star+\eps_i, \quad \text{for} \  \ i=1,\ldots
    \label{eq.lin_regression}
\end{align}
for i.i.d.\ noise variables $\eps_i \sim \mN(0,\sigma^2)$ that are also independent of the covariates $\bX_1,\ldots$ We assume that $\sigma>0$ is known. 

Let $(\bX,Y)$ be a new and independent sample with the same distribution as the training samples. Write 
\begin{align*}
    Q:=\ex[\bX \bX^\top]
\end{align*}
for the (uncentered) covariance matrix of the design vectors. Consider a possibly randomized estimator $\wh \btheta.$ By construction, $\wh \btheta$ is independent of a test point $(\bX,Y).$ Rewriting $Y-\bX^\top \wh \btheta =(Y-\bX^\top \btheta^\star) + \bX^\top (\wh \btheta-\btheta^\star),$ conditioning on $\wh \btheta$ and using tower rule, the excess risk of $\wh \btheta$ is
\begin{align}
  \ex\sbr*{(Y-\bX^\top \wh \btheta)^2 - (Y-\bX^\top \btheta^\star)^2}
  &~=~
  \ex\sbr*{2 (Y-\bX^\top \btheta^\star)  \bX^\top (\wh \btheta-\btheta^\star)+ (\wh \btheta-\btheta^\star)^\top \bX \bX^\top (\wh \btheta-\btheta^\star) }
  \notag \\
  &~=~ 
  \ex\sbr*{(\wh \btheta-\btheta^\star)^\top Q (\wh \btheta-\btheta^\star) } \notag  \\
  &~=~ 
  \Tr(QS)
  \label{eq.risk_rewrite}
\end{align}
with $S \df \ex\big[\big(\wh \btheta-\btheta^\star\big) \big(\wh \btheta-\btheta^\star\big)^\top\big].$

\cite{2023arXiv230610529C} relates gradient descent with dropout to a vector
autoregressive (VAR) process with random coefficients \cite{regis_et_al_2022}. A similar argument can be made here. Defining $\bW_k := \btheta_k - \btheta^\star,$ $G_k:= I - 2 \alpha_k \big(e^{-\bU_k}-e^{\bU_k}\big) \big(\bU_k'-\bU_k\big)^\top \bX_k  \bX_k^\top$ and $\bxi_k:= 2 \alpha_k \eps_k \bX_k^\top (\bU_k'-\bU_k) \big(e^{-\bU_k}-e^{\bU_k}\big)
  +\alpha_k \del*{ (\bX_k^\top \bU_k)^2
    - (\bX_k^\top \bU_k')^2
  } \big(e^{-\bU_k}-e^{\bU_k}\big),$ the update formula can be written in the form of a lag one VAR process \begin{align}
    \bW_k = G_k \bW_{k-1} + \bxi_k,
    \label{eq::RAR}
\end{align} with independent random coefficients $G_k$ and independent noise/innovation variables
$\bxi_k$. It can be checked that the noise is centered, that is, 
\begin{align}
    \E[\bxi_k] = 0.
    \label{eq.centered_noise_in_VAR}
\end{align}
A proof of this fact and a derivation of \eqref{eq::RAR} can be found in Appendix \ref{sec.proofs_ub}. As we can tune the learning rate $\alpha_k$ and consider the parameter $A$ from the uniform distribution as hyperparameter, it is interesting to work out the scaling of the different terms in these parameters. For small $A$, $e^{-\bU_k}-e^{\bU_k}\approx -2\bU_k$ which is of order $A.$ This means that the term $2 \alpha_k \big(e^{-\bU_k}-e^{\bU_k}\big) \big(\bU_k'-\bU_k\big)^\top \bX_k  \bX_k^\top$ in the definition of $G_k$ will be of order $\alpha_k A^2.$ The same is true for the first term of $\bxi_k$ while the second term of $\bxi_k$ scales like $O(\alpha_k A^3).$ This suggests to define 
\begin{align}
    \textit{effective learning rate} :=\alpha_k A^2.
    \label{eq.def_effective_LR}
\end{align}
This quantity will remain constant if $A$ is decreased by a factor $\gamma<1,$ and, concurrently, $\alpha_k$ is increased by a factor $\gamma^{-2}.$ Thus by making $\gamma$ and thus $A$ small, the term of the order $O(\alpha_k A^3)$ becomes negligible. While the considered optimization scheme is motivated by Hebbian learning with $A$ a small constant, treating $A$ as a hyperparameter also suggests to choose a small $A.$

The updates \eqref{eq.update_with_squared_loss} are independent of $Q=\E[\bX\bX^\top].$ As common in the literature on gradient descent, the chosen learning rate will, however, depend on the smallest eigenvalue of $Q.$

Applying \eqref{eq.risk_rewrite} to the estimator $\btheta_k$, we need to control $$S_k:=\E[(\btheta_k-\btheta^\star)(\btheta_k-\btheta^\star)^\top]=\E[\bW_k \bW_k^\top].$$ 
We can relate this to the previous iterate via $\bW_k \bW_k^\top = G_k \bW_{k-1} \bW_{k-1}^\top G_k + \bxi_k \bxi_k^\top + G_k \bW_{k-1} \bxi_k^\top + \bxi_k \bW_{k-1}^\top G_k.$ Based on this identity, we will now derive a recursive formula for $S_k$. Assume that $\bU,\bU'\sim\Unif[-A,A]^d$ are independent and also independent of $\bX.$ Set $\bD=e^{-\bU}-e^{\bU}$ and let $\bx$ be a fixed $d$-dimensional vector. The following two matrices characterize the interaction between the $d$-dimensional noise vectors $(\bU, \bU'),$
\begin{align}
    \begin{split}
  V(\bx) &~\df~
    \E\big[\bx^\top (\bU'-\bU) \bD \bD^\top (\bU'-\bU)^\top \bx\big]
  \\
  W &~\df~
    \E\big[\del*{
      (\bX^\top \bU)^2
      - (\bX^\top \bU')^2
    }^2
    \bD \bD^\top\big].
  \end{split}
  \label{eq.VZ_def}
\end{align}

\begin{lem}
\label{lem.recursion}
\begin{itemize}
    \item[(i)] the interaction terms $G_k \bW_{k-1} \bxi_k^\top$ and $\bxi_k \bW_{k-1}^\top G_k$ have mean zero,
    \item[(ii)] $\Cov(\bxi_k)=\E[\bxi_k \bxi_k^\top]=4 \alpha_k^2\sigma^2 
  \ex \sbr*{
   V(\bX)} + \alpha_k^2 W,$
   \item[(iii)] if $\mu:=-\E[U(e^{-U}-e^U)]$ with $U\sim \Unif([-A,A]),$ then,
\begin{align*}
  S_k
  &=
  \begin{aligned}[t]
    &
    \big(
    I
    - 2 \alpha_k\mu Q
  \big) S_{k-1} \big(
    I
    - 2 \alpha_k \mu Q
  \big)
  + 4 \alpha_k^2 \set*{
    \ex\sbr*{
      \bX^\top S_{k-1} \bX
      V(\bX)
    }
    - \mu^2 Q S_{k-1} Q
  } 
  \\
  & + 4 \alpha_k^2\sigma^2 
  \ex \sbr*{
   V(\bX)} + \alpha_k^2 W.
\end{aligned}
\end{align*}
and moreover, 
\begin{align}
  S_k
  &\leq 
    \|S_{k-1}\| \Big(
    I
    - 4 \alpha_k \mu Q
  + 4 \alpha_k^2
    \ex\sbr*{
      \|\bX\|^2
      V(\bX)
    }\Big) + 4 \alpha_k^2\sigma^2 
  \ex \sbr*{
   V(\bX)
  } + \alpha_k^2 W.
  \label{eq.Sk_recursion_ineq}
\end{align}
\end{itemize}
\end{lem}

In summary, $S_k$ is an affine function in $S_{k-1}.$ 

As discussed before, the last term $\alpha_k^2 W$ can be made negligible by choosing a sufficiently small parameter $A$ in the uniform distribution. The term $\mu$ is $O(A^2)$ and, under moment assumptions on $\bX,$ the respective order of the terms $\|\ex\sbr*{V(\bX)}\|$ and $\|\ex\sbr*{\|\bX\|^2V(\bX)}\|$ is $O(A^4d)$ and $O(A^4d^2).$ To decrease the norm of the error matrix $S_k$ in \eqref{eq.Sk_recursion_ineq}, the factor $I
- 4 \alpha_k \mu Q+ 4 \alpha_k^2\ex\sbr*{\|\bX\|^2V(\bX)}$ has to remain below the identity matrix. This implies that $\alpha_k A^2 Q\gtrsim \alpha_k^2A^4 d^2 I$ and motivates to choose the effective learning rate $\alpha_kA^2$ in \eqref{eq.def_effective_LR} to be bounded by $\lesssim \lambda_{\min}(Q)/d^2.$ Compared to standard choices such as $1/k,$ the learning rate in the beginning is thus small and learning makes little progress.

If $k\gtrsim d^2,$ effective learning rate $1/k$ is possible by choosing $\alpha_k A^2\asymp 1/(k\vee\lambda_{\min}(Q)d^2).$ The rate is then determined by the contributions $4 \alpha_k^2\sigma^2 \|\ex \sbr*{V(\bX)}\|=O(\sigma^2 d/k^2).$ The $1/k^2$ becomes a $1/k$ as the contributions from the previous $O(k)$ iterates add up. Up to a $\log^2(d)$-factor, this explains the convergence rate stated in \eqref{eq.Sk_operator_norm_bd} below.

To simplify the exposition, we work now with specific choices for the parameter $A$ in the uniform distribution and the learning rate $\alpha_k.$ The expressions also depend on the variance of the noise $\sigma^2$ in the linear regression model, which is assumed to be known. To obtain a convenient expressions for the upper bound, we moreover introduce
$$M_k:=\max_{i=1,\ldots,d} \ 1\vee E\big[X_i^k\big],$$ where $X_i$ denotes the $i$-th component of $\bX.$

\begin{thm}\label{thm.rate_ub}
Assume $d\geq 9.$ Let $A=\sigma/\sqrt{d}$, and choose the learning rate 
\begin{align}
    \alpha_k = \frac{11 B\log(d)}{A^2 \lambda_{\min}(Q)(Bk+d^2\log(d))} \quad \text{with} \ \ B:= 1\wedge \frac{\lambda_{\min}(Q)^2}{2904 M_4}.
    \label{eq.alpha_k_specific}
\end{align}
Then there exists a constant $C=C(M_4,\lambda_{\min}(Q)),$ such that for all $k>2d^2\log(d)/B,$ we have
\begin{align}
    \big\|\E[(\btheta_k-\btheta^\star)(\btheta_k-\btheta^\star)^\top]\big\|
    \leq \Big(\frac{2d^2\log(d)}{Bk}\Big)^{\log(d)} \|S_0\|
    + C\frac{\sigma^2 d\log^2(d)}{k}.
    \label{eq.Sk_operator_norm_bd}
\end{align}
If moreover $\Tr(Q)\leq d,$ the excess risk \eqref{eq.risk_rewrite} is bounded by
\begin{align*}
    \ex\sbr*{(\btheta_k-\btheta^\star)^\top Q (\btheta_k-\btheta^\star) }
    \leq d\Big(\frac{2d^2\log(d)}{Bk}\Big)^{\log(d)} \|S_0\|
    + C\frac{\sigma^2 d^2\log^2(d)}{k}.
\end{align*}
\end{thm}

Assuming $\Tr(Q)\leq d$ is natural and includes in particular the case that $Q=\E[\bX\bX^\top]$ is the identity matrix. 

Initializing $\btheta_0 = \bm{0}$ with the zero vector, we have $S_0=\btheta^\star (\btheta^\star)^\top$ and with Lemma \ref{lem.spect_norm_of_vec_prod}, $\|S_0\|=\|\btheta^\star\|^2.$

The rate consists of two terms. For any $\kappa>0, d>e^\kappa,$ and $k \geq 2e^{\gamma + 3\kappa} d^2\log(d)/B,$ we have 
\begin{align}
    \Big(\frac{2d^2\log(d)}{Bk}\Big)^{\log(d)}\lesssim \frac {1}{d^\gamma k^\kappa}.
    \label{eq.first_term_asymp}
\end{align} 
A proof of this inequality is given in Appendix \ref{sec.proof_ub_convergence rate}. The inequality states that for sufficiently large $k\gtrsim d^2\log(d),$ the first term is negligible and the rate is $\sigma^2d^2\log^2(d)/k.$ This rate is slower than the usual rate $\sigma^2 d/k$ if the full data are observed in the linear regression model. The $\log^2(d)$ factor seems to be an artifact of the proof and we will later derive nearly-matching lower bounds. 

Why do we loose a factor $d$ due to querying? It is tempting to assume that observed random variables from the query model and the linear regression model are equally informative. Then one can link the additional factor $d$ for the query model to the loss of information in the data: While in the linear regression model we observe in every round a $d$-dimensional covariate vector $\bX_k$ together with a real-valued response $Y_k,$ the query model observes in each round two real-valued queries. The total number of observed random variables in the query model is thus decreased by an order $O(d)$. Differently speaking the query model needs $O(d)$ more iterations to receive the same number of observed variables. The concepts in \cite{pmlr-v49-steinhardt16} might be suitable to formalize such an argument.

The result is similar in spirit to the bounds obtained in \cite{BosSH2023} for another biologically motivated learning rule called (weight-perturbed) forward gradient descent \cite{baydin2022gradients, ren2022scaling}. For a sequence of i.i.d.\ random vectors $\bzeta_1,\bzeta_2,\ldots \mN(0,I_d),$ forward gradient descent with learning rate $\alpha_k'$ is given by the update rule
\begin{equation*}
    \btheta_k=\btheta_{k-1}-\alpha_k'\big(\nabla L(\btheta_{k-1})\big)^\top \bzeta_k \bzeta_k, \quad k=1,2,\ldots
\end{equation*}
The similarity is best explained for squared loss $L.$ In this case, $\nabla L(\btheta) 
    = -2\big(Y_k-\bX_k^\top \btheta \big)\bX_k.$ If instead of an independent draw, we set $\bU_k':=-\bU_k,$ then the optimization scheme \eqref{eq.update_with_squared_loss} becomes
\begin{align*}
  \btheta_k
  &=
  \btheta_{k-1} + \alpha_k \del*{
    \big(Y_k - \bX_k^\top (\btheta_{k-1} + \bU_k)\big)^2
    - \big(Y_k - \bX_k^\top (\btheta_{k-1} - \bU_k)\big)^2
  } \big(e^{-\bU_k}-e^{\bU_k}\big)\\
  &=
  \btheta_{k-1} + 2\alpha_k \nabla L(\btheta_{k-1})^\top \bU_k\big(e^{\bU_k}-e^{-\bU_k}\big).
\end{align*}
For small $A,$ $e^{\bU_k}-e^{-\bU_k}\approx 2\bU_k,$ which means that $\btheta_k\approx \btheta_{k-1} + 4\alpha_k \nabla L(\btheta_{k-1})^\top \bU_k \bU_k.$ Whereas forward gradient descent has been proposed for normally distributed $\bzeta_k,$ the gradient-free optimization scheme here is closely related to the case where these vectors are sampled from a uniform distribution. The choice $\bU_k':=-\bU_k$ leads to less remainder terms in the analysis and less additional noise, but we find this choice less appealing to model Hebbian learning in the brain.

\section{Lower bound for adaptive queries}
\label{sec.lb_adaptive}

We now derive nearly-matching lower bounds. To make the lower bound rigorous, we need to formalize the query model. The observed data consists of $2k$ query vectors $\bv_1,\bv_1',\ldots,\bv_k,\bv_k'$ and a random vector $(Z_1,Z_1',\ldots,Z_k,Z_k')$ of length $2k.$

The unobserved/latent i.i.d.\ pairs $(\bX_\ell,Y_\ell),$ $\ell=1,\ldots,k$ are generated from the linear regression model \eqref{eq.lin_regression}. For the lower bounds we will moreover assume that the design distribution $P_\bX$ is $\mN(0,I_d)$ with $I_d$ the $d\times d$ identity matrix. This implies that $Q=E[\bX\bX^\top]=I_d.$

Let $(\mG_\ell)_{\ell\geq 0}$ be a filtration generated by exogenous randomness. We assume that for any $\ell=1,2,\ldots,$ the query vectors $\bv_\ell,\bv_\ell'$ are measurable with respect to the $\sigma$-algebra
\begin{align*}
    \mF_{\ell-1}:=\sigma\big(Z_1,Z_1',\ldots,Z_{\ell-1},Z_{\ell-1}',\bv_1, \bv_1',\ldots,\bv_{\ell-1},\bv_{\ell-1}'\big)\times \mG_{\ell-1}.
\end{align*}
This means that the query vectors might depend on past queries, past query vectors and exogenous randomness. Given $\bv_\ell,\bv_\ell',$ we observe in the $\ell$-th step the queries $(Z_\ell,Z_\ell')$ with 
\begin{align}
    Z_\ell=Y_\ell-\bX_\ell^\top \bv_\ell, \quad \quad 
    Z_\ell'=Y_\ell-\bX_\ell^\top \bv_\ell'.
    \label{eq.zlzl'}
\end{align}
Finally, an estimator $\widehat \btheta_k$ is any measurable function in $\mF_k.$ This means that an estimator is allowed to depend on all observed queries, all observed query vectors and the exogenous randomness.

To allow for two queries instead of one makes the lower bounds considerable more involved as $Z_\ell$ and $Z_\ell'$ are dependent by virtue of sharing $\bX_\ell,Y_\ell$. One query is, however, not enough to realize the zeroth-order scheme \eqref{eq.update_with_squared_loss}.

Before stating the lower bounds, we derive an equivalent representation of this model. Two statistical models $(P_{\btheta}:\btheta\in \Theta)$ and $(Q_{\btheta}:\btheta\in \Theta)$ with the same parameter space $\Theta$ are equivalent if there exist Markov kernels $M, M'$ that are independent of unknown parameters such that $Q_{\btheta}=MP_{\btheta}$ and $P_{\btheta}=M'Q_{\btheta}$ for all $\btheta\in \Theta.$ 

\begin{lem}
\label{lem.equiv_model}
Consider the query model $\{Z_1,Z_1',\ldots,Z_k,Z_k', \bv_1,\bv_1',\ldots,\bv_k,\bv_k'\}$ with \eqref{eq.zlzl'} replaced by 
\begin{align}
    Z_\ell=Y_\ell-\bX_\ell^\top \bv_\ell, \quad \quad 
    Z_\ell'=\bX_\ell^\top \bv_\ell', \quad\quad\quad\quad \ell=1,2,\ldots
    \label{eq.zlzl'_2}
\end{align}
Both models are statistically equivalent.
\end{lem}

For the lower bounds we will work in the transformed query model \eqref{eq.zlzl'_2}. If we choose the query vector $\bv_\ell'=0,$ then $Z_\ell'=0$ and this is as informative as only querying the model once.

For the model with queries \eqref{eq.zlzl'_2}, we have
\begin{align}
    \left(
    \begin{array}{c}
    Z_\ell \\
    Z_\ell' 
    \end{array}
    \right) \bigg | \, \bv_{\ell}, \bv_\ell' \sim 
    \mN
    \left(
    \left(
    \begin{array}{c}
    0 \\
    0 
    \end{array}
    \right)
    ,
    \left(
    \begin{array}{cc}
     \sigma^2+\|\btheta-\bv_\ell\|^2    &  \langle \btheta-\bv_\ell,\bv_\ell'\rangle \\
     \langle \btheta-\bv_\ell,\bv_\ell'\rangle & \|\bv_\ell'\|^2
    \end{array}
    \right)
    \right).
    \label{eq.zlzl'_distrib_2}
\end{align}
using that, we also assumed $\bX_k\sim \mathcal{N}(0,I)$ and thus $Q=I_d.$

The choice of the query vectors up to round $k,$ will be called the {\it query strategy} (up to round $k$) and is denoted by $\mV_k.$ It is important to recall that the query strategy cannot depend on $\btheta$ as it is unknown.

The data distribution of $(Z_1,Z_1',\ldots,Z_k,Z_k', \bv_1,\bv_1',\ldots,\bv_k,\bv_k')$ will be denoted by $P_{\btheta,\mV_k}$ and depends on the underlying parameter $\btheta$ and the query strategy $\mV_k.$ Since the choice of the query strategy $\mV_k$ is part of the estimation procedure, the right notion of the minimax estimation risk is 
\begin{align*}
    \inf_{\widehat \btheta, \ \mV_k} \ 
 \sup_{\btheta\in \Theta} \ \ex_{\btheta,\mV_k} \ \big[\|\widehat \btheta-\btheta\|^2\big],
\end{align*} 
where the infimum is taken over all query strategies and all estimators. The parameter space $\Theta$ will be chosen as the Euclidean ball with radius $R$, that is, $$B_R(0):=\{\btheta:\|\btheta\|\leq R\}.$$

\begin{thm}
\label{thm.lb_adapt}
If $d\geq 3$ and $k\geq d^2$, then, 
    \begin{align*}
    \inf_{\widehat \btheta, \ \mV_k} \  
 \sup_{\btheta\in B_{R}(0)} \ \ex_{\btheta,\mV_k} \ \big[\|\widehat \btheta-\btheta\|^2\big]
    \geq \frac{1}{162}\Big(1-\frac{1}{\sqrt{2}}\Big) \bigg( R^2\wedge \frac{d^2}{k}\sigma^2\bigg).
\end{align*}
\end{thm}

The statement assumes that $k\geq d^2.$ A similar condition ($k>2d^2\log(d)/B$) appears also in the corresponding upper bound in Theorem \ref{thm.rate_ub}.

Interestingly, the rate does not depend on the radius $R$ for large $k.$ This indicates that the statistical main difficulty is to recover the direction of the true regression vector.

The same upper and lower bounds hold if instead of \eqref{eq.queries_def}, we observe the squared queries $(Z_1^2,(Z_1')^2,\ldots,Z_k^2,(Z_k')^2).$ Indeed the squares are sufficient to implement the biologically inspired updating rule \eqref{eq.update_with_squared_loss}. Since $(Z_1^2,(Z_1')^2,\ldots,Z_k^2,(Z_k')^2)$ is at most as informative as $(Z_1,Z_1',\ldots,Z_k,Z_k'),$ the lower bounds remain true.

\section{Minimax risk for nonadaptive queries}
\label{sec.nonadapt}

What is the advantage to use previous information to select the next query vector? To answer this, we now consider query vectors $\bv_1,\ldots,\bv_k$ that are chosen before the data are revealed. Denote by $\mM_k$ the space of all such query strategies. For convenience, we will moreover assume throughout this section that the design distribution is standard multivariate normal, this means that 
\begin{align}
    \bX_1,\bX_2,\ldots \sim P_{\bX}=\mN(0,I_d), \ \ \text{i.i.d.}
    \label{eq.assump_design}
\end{align}
We show in this section that the minimax estimation risk with non-adaptive query strategies in $\mM_k$ and parameter space $\Theta$ the Euclidean ball $B_R(0)$ is
\begin{align}
     \inf_{\widehat \btheta , \  \mV_k \in \mM_k} \ \sup_{\btheta\in B_R(0)} \ \ex_{\btheta,\mV_k} \ \big[\|\widehat \btheta-\btheta\|^2\big] \asymp R^2\wedge \frac{d^2}{k} (R\vee \sigma)^2.
     \label{eq.minimax_non_adapt}
\end{align}

For the upper bound we construct the following estimator. If $k\leq 2d^2((\sigma^2/R^2)\vee 1),$ take $\widehat \btheta=0.$ If $k> 2d^2((\sigma^2/R^2)\vee 1),$ we have $d(k/(2d)+1)=k/2+d\leq k$ and can therefore partition the index set $\{1,\ldots,k\}$ into $d$ blocks $\mB_1,\ldots,\mB_d$ such that each block has cardinality $\geq k/(2d).$ Set $\bv_j=(R\vee \sigma)\be_s$ if $j\in B_s$ with $\be_s$ the $s$-th standard basis vector. Thanks to \eqref{eq.assump_design}, the data are then given by 
\begin{align*}
    Z_j \sim \mathcal{N}\big(0,\sigma^2+\|\btheta - (R\vee \sigma) \be_s\|^2\big)
    = \mathcal{N}\Big(0,\sigma^2+\|\btheta\|^2-2(R\vee \sigma)\btheta_s+(R\vee \sigma)^2\Big), \quad \text{if} \ \  j\in B_s
\end{align*}
and the estimator for the $s$-th component of $\btheta$ is in this case
\begin{align}
    \wh \theta_s 
    := \frac{\sigma^2+\|\btheta\|^2+(R\vee \sigma)^2-|\mB_s|^{-1}\sum_{r\in \mB_s} Z_r^2}{2(R\vee \sigma)}, \quad \text{for} \ s=1,\ldots,d,
    \label{eq.est_nonadaptive_def}
\end{align}
with $|\mB_s|$ the cardinality of the set $\mB_s.$ Whenever $r\in \mB_s,$ we have $\ex[Z_r^2]=\sigma^2+\|\btheta\|^2-2(R\vee \sigma)\btheta_s+(R\vee \sigma)^2,$ implying that $$\wh \btheta =(\wh \theta_1,\ldots, \wh \theta_d)$$ is an unbiased estimator for $\btheta.$

\begin{thm}
\label{thm.ub_non_adapt}
Assume \eqref{eq.assump_design}. For the estimator $\wh \btheta$ defined componentwise in \eqref{eq.est_nonadaptive_def}, we have
    \begin{align*}
    \sup_{\btheta\in B_{R}(0)} \ \ex_{\btheta,\mV_k} \ \big[\big\|\widehat \btheta-\btheta\big\|^2\big]
    \leq 25 \bigg( R^2\wedge \frac{d^2}{k} (R\vee \sigma)^2\bigg).
\end{align*}
\end{thm}

The estimator uses knowledge of $R$ and $\sigma^2.$ If these quantities are unknown, estimation seems unequal harder, in particular if $R$ is small.

The query vectors $\bv_j$ can be thought of as test functions or features for $\btheta.$ In the construction of the estimator, they have norm $\|\bv_j\|=(R\vee \sigma).$ Interestingly, if $\sigma>R,$ the norm exceeds $R$ and the $\bv_j$ are themselves not in the parameter space $\Theta=B_R(0).$

\begin{proof}
For $k\leq 2d^2((\sigma^2/R^2)\vee 1),$ $\wh \btheta=0$ and the result follows since $\|\wh \btheta-\btheta\|^2=\|\btheta\|^2\leq R^2.$

It remains to show that for $k> 2d^2((\sigma^2/R^2)\vee 1),$ the rate is bounded by $25\tfrac{d^2}{k}(R\vee \sigma)^2.$ The bias-variance decomposition yields 
\begin{align*}
    \ex_{\btheta,\mV_k} \big[\big\|\widehat \btheta-\btheta\big\|^2\big]=\sum_{s=1}^d \ex_{\btheta,\mV_k} \big[\big(\widehat \btheta_s-\btheta_s\big)^2\big]= 
    \sum_{s=1}^d \Bias_{\btheta,\mV_k}^2\big(\wh \btheta_s\big)
    + \sum_{s=1}^d \Var_{\btheta,\mV_k}\big(\wh \btheta_s\big).
\end{align*}
As we have already shown that $\wh \btheta$ is unbiased, it remains to bound the variances. For $\xi \sim \mN(0,a^2),$ we have $\Var(\xi^2)=\ex[\xi^4]-\ex^2[\xi^2]=3a^4-a^2=2a^2.$ Using the definition of the estimator in \eqref{eq.est_nonadaptive_def}, the independence of $Z_1,\ldots,Z_k$, that $\|\btheta\|^2\leq R^2$, that $\sigma^2+\|\btheta\|^2-2R\btheta_s+R^2\leq 5(R\vee \sigma)^2,$ and that by construction of the blocks $|\mB_s|\geq k/(2d),$ we obtain for any $j\in \mB_s,$
\begin{align*}
    \Var_{\btheta,\mV_k}\big(\wh \btheta_s\big)
    = \frac{\Var(Z_j^2)}{(2(R\vee \sigma))^2|\mB_s|} 
    = \frac{2(\sigma^2+\|\btheta\|^2-2R\btheta_s+R^2)^2}{4(R\vee \sigma)^2|\mB_s|}
    \leq \frac{50(R\vee \sigma)^4}{4(R\vee \sigma)^2|\mB_s|}
    \leq \frac{25d}{k} (R\vee \sigma)^2.
\end{align*}
Summing over $s=1,\ldots,d$ gives another factor $d$ and thus the claim follows. 
\end{proof}

We now state the corresponding lower bound.

\begin{thm}
\label{thm.non_adapt}
Assume \eqref{eq.assump_design}. If $d\geq 6,$ then for any $k=1,2,\dots$
    \begin{align*}
    \inf_{\widehat \btheta , \  \mV_k \in \mM_k} \ \sup_{\btheta\in B_{R}(0)} \ \ex_{\btheta,\mV_k} \ \big[\|\widehat \btheta-\btheta\|^2\big]
    \geq 2^{-18} \bigg( R^2\wedge \frac{d^2}{k} (R\vee \sigma)^2\bigg).
\end{align*}
\end{thm}

Together with Theorem \ref{thm.ub_non_adapt}, this shows \eqref{eq.minimax_non_adapt}. While in the adaptive setting, we had to impose the restriction $k \gtrsim d^2\log(d)$ for the upper bound and  $k\geq d^2$ for the lower bound, the derived rate in the non-adaptive setting holds for all sample sizes $k.$ 

Compared to the adaptive case, proving the lower bound in the non-adaptive case is considerably more involved. Reasons are that in this case, one additional regime occurs in the rate. To deal with this regime requires to show that whatever the choice of the query vectors is, one can find a regression vector $\btheta^\star$ in the parameter space that is far away. The fact that we allow for two queries per sample $(\bX_\ell,Y_\ell)$ makes that `far away' has to be interpreted with respect to some tube that is generated by the pair of query vectors $(\bv_\ell,\bv_\ell').$

We believe that the constant $2^{-18}$ in the lower bound can be improved significantly at the expense of a more technical proof.

The upper bound only needs one query per iteration and the lower bound is derived for two queries per iteration. This already proves that if we can only query once in every iteration, the minimax rate remains the same.

Based on the derived lower bound, we can now quantify the gap between adaptive and non-adaptive design. A natural setting is to allow that all parameters are of order one. This means that $\|\btheta\|^2=\sum_{j=1}^d \theta_j^2$ is of order $O(d)$ and motivates to choose $R=\sqrt{d}.$ 

\begin{cor}
\label{cor.gap}
Assume \eqref{eq.assump_design}. If $\sigma=1,$ $Q=I, \btheta_0=0,$ and $R=\sqrt{d}$, then for all $k\geq 2e^5d^2\log(d)/B$ (with $B$ the constant in \eqref{eq.alpha_k_specific} ) and $d\geq 8,$ the upper for the adaptive design yields the convergence rate
\begin{align*}
    \frac{d^2\log^2(d)}{k}
\end{align*}
while the minimax rate for the non-adaptive design is 
\begin{align*}
    \frac{d^3}{k}.
\end{align*}
\end{cor}

The result implies that constructing queries based on previously seen data improves the rate by a factor $d^{-1}$ (up to logarithms). The improvement will become even more pronounced if $R$ increases. Indeed, the upper bound in the adaptive query setting will remain $d^2\log^2(d)/k,$ while the minimax rate for the non-adaptive query setting becomes $d^2 R^2/k.$ We do not have a convincing heuristic argument explaining the gap in the rates. However, it is clear that the adaptive query setting can learn over time about the direction of the true $\btheta^*,$ whereas in the non-adaptive query setting one has to spread out the query vectors equally over all possible directions. This is also clearly visible in the construction of the estimator in \eqref{eq.est_nonadaptive_def}.

For the related problem of adaptive sensing, it has been found in \cite{6289365} that adaptation improves the rate by at most log-factors. In this setting there are no queries and one can choose in the $i$-th iteration a design vector $\bX_i$ based on past observations and will then observe $Y_i=\bX_i^\top \btheta^\star +\eps_i$ with independent $\eps_i\sim\mN(0,\sigma^2).$ If $s$ is the number of non-zero components of $\btheta^\star$, it is shown that the risk $\E[\|\wh \btheta_k-\btheta^\star\|^2]$ of any estimator $\wh \btheta$ based on $k$ measurements is lower bounded by $\gtrsim \sigma^2 s/k.$ On the contrary, in the non-adaptive setting with $\bX_i$ chosen i.i.d.\ and independent of previous data, the Dantzig selector $\wh\btheta_k^D$ based on $k\gtrsim s\log(d/s)$ measurements achieves 
estimation rate $\E[\|\wh \btheta_k^D-\btheta^\star\|^2]\lesssim \sigma^2 s \log(d)/k.$ Thus the gain in the convergence rate of an adaptive sampling strategy is here at most a factor $\log(d)$ in the rate. The reason why adaptation hardly improves the rate in this setting is attributed in \cite{6289365} to the difficulty to recover the support of the sparse regression vector $\btheta^\star$.

\section{Related literature}
\label{sec.related_lit}

In zeroth-order stochastic convex optimization the task is to learn a minimizer of an unknown convex function $f(\btheta) = \ex_\Xi[f(\btheta, \Xi)].$ The vanilla framework is to sequentially issue queries $\btheta_1, \btheta_2, \ldots$ and receive noisy observations $f(\btheta_k, \Xi_k)$ for i.i.d.\ unobserved $\Xi_1, \Xi_2, \ldots$. The particular case of linear regression $Y=\bX^\top \btheta^* +\epsilon$ corresponds to $\Xi = (\bX,\epsilon)$ with $d$-dimensional covariate vectors $\bX$ drawn from an unknown distribution, noise $\epsilon$, and feedback $f(\btheta,(\bX,\epsilon)) = (\bX^\top (\btheta^\star - \btheta)  + \epsilon)^2=(Y-\bX^\top \btheta)^2$. In the literature (we refer to \cite{MR3963507, MR4376588} for excellent surveys), rates for algorithms and lower bounds are organised based on three main distinctions:
\begin{itemize}
\item[-] in the \emph{optimization} literature the objective is the gap $f(\widehat \btheta_k) - \min_{\btheta} f(\btheta)$ of a proposed evaluation point $\widehat \btheta_k$, while in the \emph{bandit} literature the objective is the regret of the queries issued $\sum_{t=1}^k \del*{f(\btheta_t) - \min_{\btheta} f(\btheta)}$. The latter bounds the former by online to batch conversion, but there are interesting separations in minimax rates \cite{pmlr-v30-Shamir13}.
\item[-] in \emph{one-point feedback} the learner issues a point $\btheta$ and observes $f(\btheta, \Xi)$, where the noise $\Xi$ is i.i.d.\ between queries \cite{flaxman2004online}. In \emph{two-point feedback} the learner issues a pair of points $\btheta,\btheta'$ and observes $f(\btheta, \Xi)$ and $f(\btheta', \Xi)$ with shared noise $\Xi.$ Different rates can occur between one and two-point feedback. 
\item[-] in addition to convexity of $f$, one may get better rates by assuming that $f(\cdot)$ (or $f(\cdot, \xi)$ for every $\xi$) are Lipschitz, smooth, higher-order $\beta$-smooth \cite{bach2016highly}, and/or strongly convex.
\end{itemize}

In \cite{Duchi2015}, upper and lower bounds are derived for the gap $f(\widehat \btheta_k) - \min_{\btheta \in \Theta} f(\btheta).$ Under the imposed conditions and ignoring the dependence on the number of parameters $d$, the gap decreases in $k$ with the rate $1/\sqrt{k},$ which is slower than the $1/k$ rate obtained here. The lower bounds are obtained for linear $f(\btheta)=\E[\btheta^\top \bX]=\btheta^\top \E[\bX]$ with unknown distribution $\bX.$ In this case, the Hessian is zero and the minimizer will lie at the boundary of the parameter space $\Theta.$ Therefore the rates strongly depend on the choice of $\Theta.$

In the case of linear regression, the function $f(\btheta)=E[(Y-\bX^\top \btheta)^2]$ is a convex quadratic defined on the entire Euclidean space $\mathbb R^d$, and as such strongly convex and smooth (of infinite order) but not Lipschitz. For one-point feedback, \cite{novitskii2022improved} proves that, under suitable assumptions, the averaged iterations $\ol \btheta_k =\tfrac 1k \sum_{\ell=1}^k \btheta_\ell$ converge with rate $\tfrac{d^2}k (\tfrac k d)^{1/\beta}$ for any $\beta>0.$

Keeping the dimension $d$ fixed and letting the number of iterations tend to infinity, \cite{2021arXiv210205198J} derives a CLT for the average over all iterates (Ruppert-Polyak average). 

To complete this literature overview, we briefly mention related approaches. \cite{pmlr-v30-Shamir13} considers a zeroth-order method to learn a minimizer of the function $F(\bw)=\bw^\top A \bw +\bb^\top \bw +c$ for unknown $d\times d$ matrix $A,$ $d$-dimensional vector $\bb$ and scalar $c.$ In every iteration, one can query the function $F$ once. This is not a statistical task, as there are no data. Under suitable conditions, the considered zeroth-order method achieves the rate $d^2/k,$ with $k$ the number of iterations and it is moreover shown that this rate is optimal. Furthermore, optimal rates are also known in the case that $f$ is strongly convex and smooth \cite{agarwal2010optimal, pmlr-v30-Shamir13}.

These rates can further be contrasted to those for stochastic first-order feedback, where the feedback for a query $\theta$ is the stochastic gradient $\nabla f(\theta, \Xi)$. Here \cite{bach2013non} show a gap rate of order $\frac{\sigma^2 d + \norm{\hat\theta_0 - \theta_*}^2}{k}$ for the average iterate of SGD with constant step size. In \cite{dieuleveut2017harder} these rates are further improved with acceleration to order $\frac{\sigma^2 d}{k} + \frac{\norm{\hat\theta_0 - \theta_*}^2}{k^2}$, matching lower bounds in both contributions. The results were extended beyond the least squares setting in \cite{LSAHowFarConstantGo} and the related problem of logistic regression was analyzed in the stochastic optimization \cite{bach2016highly} and bandit settings \cite{hazan2014logistic}. For logistic regression the function $f$ is Lipschitz and higher-order smooth but not strongly convex.

\appendix

\section{Proofs for the upper bound}
\label{sec.proofs_ub}

{\it Proof of \eqref{eq::RAR}:} Expanding the squares yields
\[
  \btheta_k
  ~=~
  \btheta_{k-1} + \alpha_k \Big(
    2 (Y_k - \bX_k^\top \btheta_{k-1}) \bX_k^\top (\bU_k'-\bU_k)
    + (\bX_k^\top \bU_k)^2
    - (\bX_k^\top \bU_k')^2
  \Big) \big(e^{-\bU_k}-e^{\bU_k}\big)
\]
Setting $\bD_k:=e^{-\bU_k}-e^{\bU_k},$ the update can be rewritten as affine function in $\btheta_k$,
\begin{align*}
  \btheta_k
  =
  \Big(
    I
    - 2 \alpha_k \bD_k (\bU_k'-\bU_k)^\top \bX_k  \bX_k^\top
  \Big) \btheta_{k-1}
  + \alpha_k \Big(
    2 Y_k \bX_k^\top (\bU_k'-\bU_k)
    + (\bX_k^\top \bU_k)^2
    - (\bX_k^\top \bU_k')^2
  \Big) \bD_k,
\end{align*}
or
\begin{align*}
  \btheta_k-\btheta^\star
  ~=~
  &\Big(
    I
    - 2 \alpha_k \bD_k \big(\bU_k'-\bU_k\big)^\top \bX_k  \bX_k^\top
  \Big) (\btheta_{k-1}-\btheta^\star)
  \\
  &~+~ 2 \alpha_k (Y_k - \bX_k^\top \btheta^\star) \bX_k^\top (\bU_k'-\bU_k) \bD_k
  \\
  &~+~ \alpha_k \del*{
    (\bX_k^\top \bU_k)^2
    - (\bX_k^\top \bU_k')^2
  } \bD_k
  .
\end{align*}
Noticing that $\eps_k=Y_k - \bX_k^\top \btheta^\star,$ \eqref{eq::RAR} follows. \qed

{\it Proof of \eqref{eq.centered_noise_in_VAR}:} We show that $\E[\bxi_k]=0.$ By definition $\bxi_k:= 2 \alpha_k \eps_k \bX_k^\top (\bU_k'-\bU_k) \big(e^{-\bU_k}-e^{\bU_k}\big)
  +\alpha_k \del*{ (\bX_k^\top \bU_k)^2
    - (\bX_k^\top \bU_k')^2
  } \big(e^{-\bU_k}-e^{\bU_k}\big).$ The first term has expectation zero, since $\E[\eps_k]=0$ and $\eps_k$ is independent of all the other variables. Since $\bU_k$ and $-\bU_k$ have the same distribution, $\E[e^{-\bU_k}-e^{\bU_k}]=0$ and $\E[(\bX_k^\top \bU_k)^2\big(e^{-\bU_k}-e^{\bU_k}\big)]=0.$ Using the independence of $\bU_k,\bU_k',$ 
  \begin{align*}
      \E[\bxi_k]
      = \E\big[\alpha_k \del*{ (\bX_k^\top \bU_k)^2
    - (\bX_k^\top \bU_k')^2
  } \big(e^{-\bU_k}-e^{\bU_k}\big)\big]
  =0.
  \end{align*}

\subsection{Expectations with respect to the uniform distribution}
\label{sec.expectation_uniform}

The moments of the noise $U$ play an important role in the analysis.

\begin{defi} For natural numbers $r,q \ge 0$, we abbreviate
\begin{align*}
  c_{r,q}
  ~:=~
  \ex \big[U^r(e^{-U}-e^U)^q\big] \quad \text{with} \ \  U\sim \Unif[-A,A].
\end{align*}
\end{defi}

We have $c_{0,0} ~=~ 1$ and  $c_{2,0} ~=~ A^2/3.$ By definition, $-U(e^{-U}-e^U)$ is a non-negative random variable and therefore $-c_{1,1}>0.$ More specifically,

\begin{lem}
\label{lem.mu_expression}
If $U\sim \Unif([-A,A])$ and $A\leq 1,$ then 
\begin{align}
    \frac{A^2}{11}\leq -c_{1,1}\leq A^2 \quad \text{and} \ \ c_{0,2}\leq - 3c_{1,1}.
    \label{eq.ublb_mu}
\end{align}
\end{lem}

\begin{proof}
Integration by parts gives
\begin{align*}
	-\E[U(e^{-U}-e^U)]
	&= 
	-\frac{1}{2A} \int_{-A}^A u\big(e^{-u}-e^u\big) \, du
	= \frac{1}{A}\int_{-A}^A u e^u \, du
	= \frac{1}{A} ue^u \, \Big|_{-A}^A-\frac 1A\int_{-A}^A e^u \, du\\
	&= e^A+e^{-A}-\frac{e^A-e^{-A}}{A}.
\end{align*}
For the second part observe that for $A\geq 0,$ third order Taylor expansion gives $e^A\geq 1+A+A^2/2$ and $e^{-A}\leq 1-A+A^2/2.$ From the expression above and using that $A\leq 1,$
\begin{align*}
	-\E[U(e^{-U}-e^U)]
	&= e^A+e^{-A}-\frac{e^A-e^{-A}}{A} \\
&= -e^A \Big(\frac 1A-1\Big)+e^{-A}\Big(1+\frac 1A\Big) \\
&\leq -\Big( 1+A+\frac{A^2}2\Big) \frac{1-A}{A}
+ \Big(1-A+\frac{A^2}2\Big) \frac{A+1}{A} \\
&= -(1+A) \frac{1-A}{A} - \frac{A-A^2}2 +(1-A) \frac{A+1}{A}
+ \frac{A^2+A}2 \\
&=A^2.
\end{align*}
By third order Taylor expansion and $A\leq 1,$ we find $e^A\leq 1+ A+A^2/2+e^A A^3/6\leq 1+ A+A^2/2+e A^3/6$ and $e^{-A}\geq 1 -A +A^2/2-e A^3/6.$ Basically following the same steps as for the upper bound of $-\E[U(e^{-U}-e^U)],$
\begin{align*}
	-\E[U(e^{-U}-e^U)]
&= -e^A \Big(\frac 1A-1\Big)+e^{-A}\Big(1+\frac 1A\Big) \\
&\geq -\Big( 1+A+\frac{A^2}2+e \frac{A^3}{6}\Big) \frac{1-A}{A}
+ \Big(1-A+\frac{A^2}2-e \frac{A^3}{6}\Big) \frac{A+1}{A} \\
&=A^2\Big(1-\frac{e}{3}\Big).
\end{align*}
Since $1-e/3\geq 0.093\geq 1/11,$ this completes the proof for \eqref{eq.ublb_mu}.

To prove $c_{0,2}\leq -3c_{1,1},$ we use that $e^x-e^{-x}=2\sum_{\ell \ \text{odd}} x^\ell/\ell!.$ Since for odd $\ell,$ $x^\ell (e^x-e^{-x})\geq 0,$ we have for $|x|\leq 1$ that $x^\ell (e^x-e^{-x})\leq x(e^x-e^{-x}).$ With $e\leq 3$, we obtain for $x\leq 1,$ 
\begin{align*}
    \big(e^x-e^{-x}\big)^2 
    &= 2\sum_{\ell \ \text{odd}} \frac{x^\ell}{\ell!} \big(e^x-e^{-x}\big)
    \leq 2\sum_{\ell \ \text{odd}} \frac{1}{\ell!} x\big(e^x-e^{-x}\big) \\
    &\leq 2\Big(e-1-\frac 1{2!} \Big) x\big(e^x-e^{-x}\big)
    \leq 3 x\big(e^x-e^{-x}\big).
\end{align*}
Thus $c_{0,2}=\E[(e^U-e^{-U})^2]\leq 3\E[U(e^U-e^{-U})]=-3c_{1,1}.$
\end{proof}

Using the previous lemma,
\begin{align}
	c_{r,2}\leq A^r c_{0,2}\leq 3A^{r+2}.
	\label{eq.cr2_bd}
\end{align}
For odd $r$, $c_{r,2}=0$ and this inequality makes only sense if $r$ is even.

\subsection{Moments of the noise contributions in the zeroth-order scheme}

Recall that $\|\cdot\|$ denotes the spectral norm for matrices and the Euclidean norm for vectors.

\begin{lem}
\label{lem.spect_norm_of_vec_prod}
Let $\ba,\bb$ be column vectors of the same length. Then $\|\ba \bb^\top\| = \|\ba\|\|\bb\|$.
\end{lem}

\begin{proof}
The matrix $\bb \bb^\top$ is of rank one with non-zero eigenvalue $\bb^\top \bb$ and corresponding eigenvector $\|\bb\|^{-1}\bb.$ Thus,
$\|\ba \bb^\top\|^2=\lambda_{\max}(\bb \ba^\top \ba \bb^\top)=\ba^\top \ba \lambda_{\max}(\bb \bb^\top)=\ba^\top \ba \bb^\top \bb.$ The result follows by taking square roots.
\end{proof}

We now derive closed-form expressions and bounds for the two expected values in \eqref{eq.VZ_def}. Recall that $c_{r,q}~=~\ex \big[U^r(e^{-U}-e^U)^q\big].$ The $k$-th power $\bx^k$ of a vector $\bx=(x_1,\ldots,x_d)$ is understood componentwise, e.g., $\norm{\bx^2}^2=\sum_i x_i^4.$
\begin{lem}
  Let $I$ denote the $d\times d$ identity matrix. We have
\begin{align*}
  V(\bx) &~=~(c_{2,2} - c_{2,0} c_{0,2} - 2 c_{1,1}^2) \diag(\bx^2)
  + 2 c_{1,1}^2 \bx \bx^\top
  + 2c_{2,0} c_{0,2} \norm*{\bx}^2 I.
  \intertext{and}
  W
  ~=~
  \E\Big[&(c_{4,2}
  - 6 c_{2,0} c_{2,2}
  + 6 c_{2,0}^2 c_{0,2}
  - c_{4,0} c_{0,2}
  - 8 c_{3,1} c_{1,1}
  + 24 c_{2,0} c_{1,1}^2
  ) \diag(\bX^4)
    \\
    & ~+~ (
    4 c_{2,0} c_{2,2}
    - 4 c_{2,0}^2 c_{0,2}
    - 8 c_{2,0} c_{1,1}^2
    ) \norm{\bX}^2 \diag(\bX^2)
    \\
    & ~+~ (
    4 c_{3,1} c_{1,1}
    - 12 c_{2,0} c_{1,1}^2
    ) \big(\bX^3 \bX^\top + \bX (\bX^3)^\top\big)
  \\
  & ~+~ 8 c_{2,0} c_{1,1}^2 \norm{\bX}^2 \bX \bX^\top
  \\
  & ~+~ (2 c_{4,0} c_{0,2} - 6 c_{2,0}^2 c_{0,2}) \norm{\bX^2}^2 I
  \\
  & ~+~ 4 c_{2,0}^2 c_{0,2} \norm{\bX}^4 I
  \Big].
\end{align*}
Moreover for $\bX=(X_1,\ldots,X_d)^\top,$
\begin{align}
    \ex[V(\bX)]
    &\leq 12A^4d \max_{i=1,\ldots,d} \ex[X_i^2]I, \label{eq.V(x)_bd} \\
    \ex\big[\|\bX\|^2 V(\bX)\big]
    &\leq12A^4d^2 \max_{i=1,\ldots,d} \ex[X_i^4] I, \label{eq.x^2V(x)_bd} \\
    W
    &\leq 107 A^6 d^2 \max_{i=1,\ldots,d} 1\vee \ex[X_i^4] I. \label{eq.Z(x)_bd}
\end{align}
\end{lem}

While the closed-form expression of $W$ depends on fourth power of $\bX,$ it is convenient to relate this to a sixth power in \eqref{eq.Z(x)_bd}.

\begin{proof}
  We first prove the formula for $W.$ Let $\bx=(x_1,\ldots,x_d)$ be fixed. Using that $\ex [(\bx^\top \bU')^2]=\ex [(\bU')^\top \bx \bx^\top \bU']=c_{2,0}\tr(\bx \bx^\top)=c_{2,0}\bx^\top \bx,$
  \begin{align*}
  \ex [(\bx^\top \bU')^4]
  &=\ex \Big[\sum_{i,j,k,\ell} U_i'U_j'U_k'U_\ell' x_ix_jx_kx_\ell \Big] \\
  &=c_{4,0}\sum_i x_i^4+3c_{2,0}^2 \sum_{i\neq j} x_i^2 x_j^2 \\
  &=(c_{4,0}-3c_{2,0}^2)\|\bx^2\|^2+3c_{2,0}^2\|\bx\|^4,    
  \end{align*}
  $\ex [\bD\bD^\top]=c_{0,2}I,$ and the independence of $\bU$ and $\bU'$ gives
\begin{align*}
  &\ex \sbr*{
    \del*{
      (\bx^\top \bU)^2
      - (\bx^\top \bU')^2
    }^2
    \bD \bD^\top
  }
  \\
  &~=~
  \ex \sbr*{
    \del*{
      (\bx^\top \bU)^4
      - 2 c_{2,0} (\bx^\top \bU)^2 \bx^\top \bx
      + \del*{c_{4,0} - 3 c_{2,0}^2} \norm{\bx^2}^2
      + 3 c_{2,0}^2 \norm{\bx}^4
    }
    \bD \bD^\top
  }
  \\
  &~=~
  \ex \sbr*{
      (\bx^\top \bU)^4
    \bD \bD^\top
  }
  -
  2 c_{2,0}  \bx^\top \bx
  \ex \sbr*{
    (\bx^\top \bU)^2
    \bD \bD^\top
  }+
  c_{0,2}
  \del*{
    \del*{c_{4,0} - 3 c_{2,0}^2} \norm{\bx^2}^2
    + 3 c_{2,0}^2 \norm{\bx}^4
  }
  I
  .
\end{align*}
Next we simplify these two remaining expectations. The $(i,j)$-th off-diagonal entry of the matrix $\ex \sbr*{(\bx^\top \bU)^2 \bD \bD^\top}$ is $\ex[\sum_{\ell,k} U_\ell U_k x_\ell x_k D_i D_j].$ The summands are non-zero if either $(\ell,k)=(i,j)$ or $(\ell,k)=(j,i).$ In both cases we get the contribution $c_{1,1}^2 x_ix_j,$ such that for $i\neq j,$  $\ex[\sum_{\ell,k} U_\ell U_k x_\ell x_k D_i D_j]=2c_{1,1}^2 x_ix_j.$ For the $i$-th diagonal entry of $\ex \sbr*{(\bx^\top \bU)^2 \bD \bD^\top}$, we obtain $\ex[\sum_{\ell,k} U_\ell U_k x_\ell x_k D_i^2]=\ex[\sum_{\ell} U_\ell^2 x_\ell^2 D_i^2]=c_{2,2} x_i^2+c_{2,0}c_{0,2}\sum_{\ell \neq i} x_\ell^2=(c_{2,2}-c_{2,0}c_{0,2})x_i^2+c_{2,0}c_{0,2} \|\bx\|^2.$ Combining these formulas yields
\begin{align}
  \begin{split}
  &\ex \sbr*{
    (\bx^\top \bU)^2
    \bD \bD^\top
  } =~
  (c_{2,2} - c_{2,0} c_{0,2} - 2 c_{1,1}^2) \diag(\bx^2)
  + 2 c_{1,1}^2 \bx \bx^\top
  + c_{2,0} c_{0,2} \norm*{\bx}_2^2 I
  \end{split}
  \label{eq.U2DD}
\end{align}
and
\begin{align*}
  \ex \sbr*{
      (\bx^\top \bU)^4
    \bD \bD^\top
  }
  &~=~
  \ex \sbr*{
    \sum_{i,j,\ell,m}
    x_i U_i x_j U_j  x_\ell U_\ell x_m U_m
    \bD \bD^\top
  }
  \\
  &~=~
  \begin{aligned}[t]
    & (c_{4,2} - 6 c_{2,0} c_{2,2} + 6 c_{2,0}^2 c_{0,2} - c_{4,0} c_{0,2} - 8 c_{3,1} c_{1,1} + 24 c_{2,0} c_{1,1}^2) \diag(\bx^4)
    \\
    & + (6 c_{2,0} c_{2,2} -6 c_{2,0}^2 c_{0,2} - 12 c_{2,0} c_{1,1}^2) \norm{\bx}^2 \diag(\bx^2)
    \\
    & + (c_{4,0} c_{0,2} - 3 c_{2,0}^2 c_{0,2}) \norm{\bx^2}^2 I
    \\
    & + 3 c_{2,0}^2 c_{0,2} \norm{\bx}^4 I
    \\
    & + 12 c_{2,0} c_{1,1}^2 \norm{\bx}^2 \bx \bx^\top
    \\
    & + (4 c_{3,1} c_{1,1} - 12 c_{2,0} c_{1,1}^2) \big(\bx^3 \bx^\top + \bx (\bx^3)^\top\big).
  \end{aligned}
\end{align*}
This is because on a diagonal entry $(p,p),$ we have $c_{4,2} x_p^4$, as well as $6 c_{2,0} c_{2,2} x_p^2 \sum_{q \neq p} x_q^2$ and $3 c_{2,0}^2 c_{0,2} (\sum_{q \neq p} x_q^2)^2$ and $c_{4,0} c_{0,2} \sum_{q \neq p} x_q^4$. On an off-diagonal entry $(p,q),$ we have $12 c_{2,0} c_{1,1}^2 x_p x_q \sum_{r \notin \set{p,q}} x_r^2$ as well as $4 c_{3,1} c_{1,1} (x_p^3 x_q + x_p x_q^3)$. Combining/grouping terms, replacing $\bx$ by the random vector $\bX$ and taking expectation with respect to $\bX$ yields the formula for $W$.

For $V(\bx),$ observe that the cross-terms are zero and by conditioning first on $\bU'$ and using that $\ex[\bD\bD^\top] =c_{0,2} I,$ we find
\begin{align*}
    V(\bx) &=
  \ex\big[
    \bx^\top (\bU'-\bU) \bD \bD^\top (\bU'-\bU)^\top \bx \big] \\
    &= \ex \big[
    \bx^\top \bU' \bD \bD^\top (\bU')^\top \bx\big]
    +\ex \big[
    \bx^\top \bU \bD \bD^\top \bU^\top \bx\big] \\
    &= c_{0,2} \bx^\top \ex \sbr*{\bU' (\bU')^\top}\bx
    +\ex \big[
    \bx^\top \bU \bD \bD^\top \bU^\top \bx\big] \\
    &= c_{0,2} c_{2,0} \|\bx\|^2 I
    +\ex \big[
    \bx^\top \bU \bD \bD^\top \bU^\top \bx\big].
\end{align*}
Combined with \eqref{eq.U2DD}, the formula $V(\bx) = (c_{2,2} - c_{2,0} c_{0,2} - 2 c_{1,1}^2) \diag(\bx^2)
  + 2 c_{1,1}^2 \bx \bx^\top
  + 2c_{2,0} c_{0,2} \norm*{\bx}^2 I$ follows.

For the bounds on the expectations, we use the bounds on the moments $c_{r,q}$ derived in Section \ref{sec.expectation_uniform}. In particular $c_{r,q}\geq 0$ whenever $r$ and $q$ are even. The matrix $\bx\bx^\top$ is positive semi-definite. By \eqref{eq.ublb_mu}, $c_{1,1}\leq A^2$ and $c_{0,2}\leq 3A^2.$ Thus, $V(\bx) \leq 3A^4\diag(\bx^2)+ 2A^4 \bx \bx^\top + 6A^4 \norm*{\bx}^2 I.$ Moreover, by Lemma \ref{lem.spect_norm_of_vec_prod}, the largest eigenvalue of $\bx\bx^\top$ is $\|\bx\|^2$ and therefore $\bx\bx^\top \leq \|\bx\|^2 I.$ Since also $\diag(\bx^2)\leq \|\bx\|^2 I,$ we obtain $V(\bx) \leq 12A^4\|\bx\|^2 I$ and thus for $\bX=(X_1,\ldots,X_d)^\top,$  $\ex[V(\bX)]\leq 12A^4d \max_{i=1,\ldots,d} \ex[X_i^2] I,$ proving \eqref{eq.V(x)_bd}. Since $X_i^2X_j^2\leq X_i^4/2+X_j^4/2,$ we can derive moreover $\ex[\|\bX\|^2 V(\bX)]\leq 12A^4d^2 \max_{i,j=1,\ldots,d} \ex[X_i^2X_j^2]\leq12A^4d^2 \max_{i=1,\ldots,d} \ex[X_i^4],$ proving  \eqref{eq.x^2V(x)_bd}.

We finally derive \eqref{eq.Z(x)_bd}. Since $2\|\bx^3 \bx^\top\|=2\|\bx^3\|\|\bx\|,$ all eigenvalues of $\bx^3 \bx^\top + \bx (\bx^3)^\top$ lie between $-2\|\bx^3\|\|\bx\|$ and $2\|\bx^3\|\|\bx\|.$ Therefore, for any real number $a,$ we have $a(\bx^3 \bx^\top + \bx (\bx^3)^\top)\leq 2|a| \|\bx^3\|\|\bx\|.$ By \eqref{eq.cr2_bd}, $c_{r,2}\leq A^r c_{0,2}\leq 3A^{r+2}.$ Moreover $0\leq c_{3,1} c_{1,1}\leq A^2 c_{1,1}^2\leq A^6$ and $0\leq c_{2,0} c_{1,1}^2\leq A^6.$ Thus,
\begin{align*}
  W
  \leq 
  &\E\Big[(3+18+24) A^6 \diag(\bX^4)
    +
    12A^6 \norm{\bX}^2 \diag(\bX^2)
    +12A^6
    2\|\bX^3\|\|\bX\| I
  +8A^6\norm{\bX}^4 I
  \\
  &+6A^6\norm{\bX^2}^2 I+12A^6\norm{\bX}^4 I\Big]
  .
\end{align*}
Since $\diag(\bX^4)\leq \|\bX\|^4 I,$ $\diag(\bX^2)\leq \|\bX\|^2 I,$ and $\|\bX^2\|^2\leq \|\bX\|^4,$ we obtain $W\leq 83 A^6 \E[\norm{\bX}^4]I+24 A^6 \E[\|\bX^3\|\|\bX\|]I.$

Because of the elementary inequality $ab\leq \tfrac 12 a^2+\tfrac 12 b^2$,
\begin{align*}
	\ex\big[\norm{\bX}^4\big]
	=\sum_{i,j=1}^d\ex[X_i^2X_j^2] \leq \frac 12 \sum_{i,j=1}^d\ex[X_i^4]+\ex[X_j^4]\leq d^2\max_{i=1,\ldots,d} \ex[X_i^4].
  \end{align*}
Similarly, $\norm{\bX^3} \norm{\bX}~\le~
      \sqrt{d \max_i X_i^6}
      \sqrt{d \max_i X_i^2}
      ~=~
      d \max_i X_i^4
      ~\le~
      d \sum_i X_i^4$ and hence
\[
      \ex\sbr*{\norm{\bX^3} \norm{\bX}}
      ~\le~
      d  \sum_i \ex[ X_i^4]
      ~\le~
      d^2 \max_i \ex[ X_i^4].
    \]      
Combined with $W\leq 83 A^6 \E[\norm{\bX}^4]I+24 A^6 \E[\|\bX^3\|\|\bX\|]I,$ we finally obtain $$W\leq 107 A^6 d^2 \max_{i=1,\ldots,d} 1\vee \ex[X_i^4] I.$$
\end{proof}

\subsection{Proof of Lemma \ref{lem.recursion}}

Recall that $\bD_k:=e^{-\bU_k}-e^{\bU_k}.$ The rewritten update equation \eqref{eq::RAR} can be decomposed into the following three parts
\begin{align} \label{eq:parts}
  \btheta_k-\btheta^\star
  &~=~
  \bA + \bB + \bC,
  \intertext{where}
  \notag
  \bA
  &~\df~ \del*{
    I
    - 2 \alpha_k \bD_k (\bU_k'-\bU_k)^\top \bX_k  \bX_k^\top
  } (\btheta_{k-1}-\btheta^\star),
  \\
  \notag
  \bB &~\df~ 2 \alpha_k \eps_k \bX_k^\top (\bU_k'-\bU_k) \bD_k,
  \\
  \notag
  \bC & ~\df~ \alpha_k \del*{
    (\bX_k^\top \bU_k)^2
    - (\bX_k^\top \bU_k')^2
  } \bD_k.
\end{align}
In the following we often use that $(\bX_k,Y_k)$ has the same distribution as $(\bX,Y).$ This means that if all randomness in an expectation is only due to $(\bX_k,Y_k)$, these variables can be replaced by $(\bX,Y).$ Recall moreover that odd powers of $\bU$ or $\bU'$ disappear and that we have defined $\mu=-c_{1,1}.$ Thus,
\begin{align}
\begin{split}
  \ex [\bA \bA^\top]
  &~=~
  S_{k-1}
  -
  2 \mu \alpha_k (Q S_{k-1} + S_{k-1} Q)
  +
  4 \alpha_k^2 \ex\sbr*{
    \bX_k^\top S_{k-1} \bX_k
    V(\bX_k)
    } \\
  &~=~ \del*{
    I
    - 2 \alpha_k \mu Q
  } S_{k-1} \del*{
    I
    - 2 \alpha_k \mu Q
  }
  + 4 \alpha_k^2 \set*{
    \ex\sbr*{
      \bX^\top S_{k-1} \bX
      V(\bX)
    }
    - \mu^2 Q S_{k-1} Q
  }.  
  \end{split}
  \label{eq.423762}
\end{align}
Since for any $d\times d$ matrix $T$, we have $T^\top S_{k-1} T\leq \|S_{k-1}\| T^\top T,$ we also have 
\begin{align}
  \ex [\bA \bA^\top]
  &~\leq~ \|S_{k-1}\| \del*{
    I
    - 4 \alpha_k \mu Q
  + 4 \alpha_k^2
    \ex\sbr*{
      \|\bX\|^2
      V(\bX)
    }}.  
    \label{eq.423762_2}
\end{align}
Now by counting powers of $\bU_k, \bU_k',\bD_k$, we find that
\[
  \ex[\bA \bC^\top]
  ~=~
  0.
\]
Applying tower rule by conditioning on all randomness except $\eps_k\sim \mN(0,\sigma^2),$ it follows that
\begin{align*}
    \ex[\bA \bB^\top]
  &~=~
  \ex\big[\ex[\bA \bB^\top \, | \, \bU_k,\bU_k',\bX_k,\btheta_k]\big]
  = \ex\Big[\bA 2\alpha_k \bD_k^\top (\bU_k'-\bU_k)^\top \bX_k\ex[\eps_k] \Big]
  = 0, \\
  \ex[\bB \bC^\top]
  &~=~
  0, \\
   \ex[\bB \bB^\top]
  &~=~
  4 \alpha_k^2
  \ex \sbr*{
         \eps_k^2 V(\bX)} = 4 \alpha_k^2\sigma^2 
  \ex \sbr*{V(\bX)},
\end{align*}
and finally
\begin{align*}
  \ex[\bC \bC^\top]
  &~=~
  \alpha_k^2 W,
\end{align*}
which is some combination of fourth powers of $\bX$. 

Now (i) follows due to $\E[G_k \bW_k \bxi_k^\top]=\E[\bA(\bB+\bC)^\top]=0.$ Statement (ii) follows since by \eqref{eq.centered_noise_in_VAR}, $\E[\bxi_k]=0$ such that $\Cov(\bxi_k)=\E[\bxi_k\bxi_k^\top]=\E[(\bB+\bC)(\bB+\bC)^\top]=\E[\bB\bB^\top]+\E[\bC\bC^\top]=4 \alpha_k^2\sigma^2 
  \ex \sbr*{V(\bX)}+\alpha_k^2 W.$ To see (iii), observe that
\begin{align*}
  S_k
  &=
  \ex[(\bA+\bB+\bC)(\bA+\bB+\bC)^\top]
  \\
  &=
  \begin{aligned}[t]
    &\ex\big[\bA \bA^\top\big] + \ex\big[\bA \bB^\top\big] + \ex\big[\bB \bA^\top\big] + \ex\big[\bB \bB^\top\big] + \ex\big[\bC \bC^\top\big]
\end{aligned}
  \\
  &=
  \begin{aligned}[t]
    &
    \big(
    I
    - 2 \alpha_k \mu Q
  \big) S_{k-1} \big(
    I
    - 2 \alpha_k \mu Q
  \big)
  + 4 \alpha_k^2 \Big(
    \ex\big[
      \bX^\top S_{k-1} \bX
      V(\bX)
    \big]
    - \mu^2 Q S_{k-1} Q
  \Big)  \\
  & + 4 \alpha_k^2\sigma^2\E[V(\bX)] + \alpha_k^2 W.
\end{aligned}
\end{align*}
Applying \eqref{eq.423762_2} instead of \eqref{eq.423762} yields the asserted inequality. \qed

\subsection{Proof of Theorem \ref{thm.rate_ub}}
\label{sec.proof_ub_convergence rate}

We now combine the recursive formula with the bounds obtained in \eqref{eq.V(x)_bd}, \eqref{eq.x^2V(x)_bd}, and \eqref{eq.Z(x)_bd}. Recall that by definition $\mu=-c_{1,1}$ and by Lemma \ref{lem.mu_expression}, $\mu\geq A^2/11$  which is the same as $-\mu\leq -A^2/11.$ Setting $M_k:=\max_{i=1,\ldots,d} 1\vee E[X_i^k],$ we obtain
\begin{align*}
  \|S_k\|
  &\leq 
    \|S_{k-1}\| \Big(
    1
    - \frac{4}{11} \alpha_k A^2 \lambda_{\min}(Q)
  + 48 \alpha_k^2 A^4d^2 M_4\Big) + 48 \alpha_k^2\sigma^2 A^4d M_2 + \alpha_k^2 107 A^6 d^2 M_4.
\end{align*}
By assumption $A= \sigma/\sqrt{d}.$ If 
\begin{align}
    96 \alpha_k^2 A^4d^2 M_4 \leq \frac{4}{11} \alpha_k A^2 \lambda_{\min}(Q)\leq 1,    
    \label{eq.alpha_k_condition}
\end{align}
then, the previous inequality yields the simpler recursion inequality
\begin{align*}
  \|S_k\|
  &\leq 
    \|S_{k-1}\| \Big(
    1
    - \frac{2}{11} \alpha_k A^2 \lambda_{\min}(Q)\Big) + (48M_2+107M_4) \alpha_k^2\sigma^2 A^4d \\
    &= \|S_{k-1}\| \big(
    1
    - \beta_k\big) + C_* \beta_k^2 \sigma^2 d,
\end{align*}
where we defined
\begin{align}
    \beta_k:=\frac 2{11} \alpha_k A^2 \lambda_{\min}(Q) \quad 
    \text{and} \ \  C_*:= (48M_2+107M_4)\Big(\frac{11}{2\lambda_{\min}(Q)}\Big)^2.
\end{align}
Induction with respect to $k$ gives 
\begin{align}
    \|S_k\|
    \leq \|S_0\| \prod_{\ell=1}^k 
    \big(
    1
    - \beta_\ell\big)
    + C_*\sum_{j=1}^k\bigg(\prod_{\ell=j+1}^k 
    \big(
    1
    - \beta_\ell \big)\bigg) \beta_j^2 \sigma^2d,
    \label{eq.Sk_spectral_ineq1}
\end{align}
where $\prod_{\ell=k+1}^{k}(1-\beta_{\ell}):=1.$ Rewritten in terms of $\beta_k,$ the assumed inequality \eqref{eq.alpha_k_condition} is  $$96\Big(\frac{11\beta_k^2}{2\lambda_{\min}(Q)}\Big)^2d^2M_4\leq 2\beta_k\leq 1.$$ The choice of the learning rate for $\alpha_k$ in \eqref{eq.alpha_k_specific} results in 
\begin{align}
    \beta_k= \frac 2{11} \alpha_k A^2 \lambda_{\min}(Q) = \frac{2B\log(d)}{Bk+d^2\log(d)} \quad \text{with} \ \ B= 1\wedge \frac{\lambda_{\min}(Q)^2}{2904M_4}.
\end{align}
For all $d\geq 1,$ we have $\beta_k\leq 2B/d^2.$ Using that $96(11/2)^2=2904$, one can now check that for all $d\geq 2,$ the assumed inequalities for $\beta_k$ hold and afortiori thus also \eqref{eq.alpha_k_condition}.

Let $k^*$ be the smallest integer such that $Bk^*\geq d^2\log(d).$ Since $d\geq 2,$ we must have 
\begin{align}
	k^*\leq \frac{1+d^2\log(d)}{B}\leq \frac{2d^2\log(d)}{B}.
	\label{eq.kstar_ub}
\end{align}
For all $k\geq k^*,$ we have $\beta_k\geq \log(d)/k.$ By bounding $1-\beta_\ell\leq 1$ for all $\ell < k^*,$ and using for $k\geq j+1\geq k^*,$ $\sum_{\ell=j+1}^k \tfrac{1}{\ell}\geq \sum_{\ell=j+1}^k \int_{\ell}^{\ell+1} \tfrac{1}{\ell} \, du = \int_{j+1}^{k+1} \frac{1}{u} \, du = \log(k+1)-\log(j+1)=\log((k+1)/(j+1)),$ we obtain for all $k\geq k^*,$
\begin{align*}
    \prod_{\ell=j+1}^{k} 
    \big(
    1
    - \beta_\ell \big)
    \leq \exp\Big(-\log(d)\sum_{\ell=k^*\vee (j+1)}^k \frac{1}{\ell}\Big)
    \leq \Big(\frac{k^*\vee (j+1)}{k+1}\Big)^{\log(d)}.
\end{align*}
Combined with $\beta_k\leq 2B/d^2$ for all $k,$ $\beta_j\leq \log(d)/j\leq 2\log(d)/(j+1)$ for all $j\geq k^*,$ and \eqref{eq.Sk_spectral_ineq1}, we obtain for all $k>k^*,$
\begin{align*}
	\|S_k\|
    &\leq \Big(\frac{k^*}{k}\Big)^{\log(d)} \|S_0\|
    + C_*\sum_{j=1}^{k^*-1} \Big(\frac{k^*}{k}\Big)^{\log(d)}
    \Big( \frac{2B}{d^2}\Big)^2 \sigma^2d
    +C_*
    \sum_{j=k^*}^k \Big(\frac{j+1}{k+1}\Big)^{\log(d)} \Big(\frac{2\log(d)}{j+1}\Big)^2 \sigma^2d \\
    &\leq \Big(\frac{k^*}{k}\Big)^{\log(d)} \|S_0\|
    + C_*\frac 1k \Big(\frac{k^*}{k}\Big)^{\log(d)-1} (k^*)^2 
    \Big( \frac{2B}{d^2}\Big)^2 \sigma^2d
    +
    C_*\sum_{j=k^*}^{k} \Big(\frac{j+1}{k+1}\Big)^{\log(d)} \Big(\frac{2\log(d)}{j+1}\Big)^2 \sigma^2d.
\end{align*}
Since $d\geq 9$, we have $\log(d)-2>0$ and therefore $k^*/k\leq 1.$ Using also  \eqref{eq.kstar_ub} gives
\begin{align*}
    \frac 1k \Big(\frac{k^*}{k}\Big)^{\log(d)-1} (k^*)^2 
    \Big( \frac{2B}{d^2}\Big)^2 \sigma^2d
    \leq \frac 1k \Big(\frac{2d^2\log(d)}{B}\Big)^2 \Big( \frac{2B}{d^2}\Big)^2 \sigma^2d
    \leq \frac{16\sigma^2 d\log^2(d)}{k}.
\end{align*}
Applying $\log(d)-2>0$ again and moreover, $(j+1)/(k+1)\leq 1$ for all $j\leq k,$ we find
\begin{align*}
    \sum_{j=k^*}^k \Big(\frac{j+1}{k+1}\Big)^{\log(d)} \frac1{(j+1)^2}=\frac{1}{(k+1)^2} \sum_{j=k^*}^k \Big(\frac{j+1}{k+1}\Big)^{\log(d)-2}\leq \frac 1k.
\end{align*}
Combined, the last three displayed inequalities combined yield
\begin{align*}
	\|S_k\|
    &\leq \Big(\frac{k^*}{k}\Big)^{\log(d)} \|S_0\|
    + C_*\frac{20\sigma^2 d\log^2(d)}{k}.
\end{align*}
The first claim follows now with $C=20C_*.$

For two positive semi-definite matrices $D=D_1^\top D_1$ and $E,$ we know that $D_1 E D_1^\top$ is positive semi-definite and bounded by $\|E\| D_1 D_1^\top.$ Therefore, $\Tr(DE)=\Tr(D_1 E D_1^\top)\leq \Tr(\|E\| D_1 D_1^\top)=\|E\|\Tr(D_1^\top D_1)=\|E\|\Tr(D).$ Applying this to $(D,E)=(Q,S_k)$ and using that $$\ex\sbr*{(\btheta_k-\btheta^\star)^\top Q (\btheta_k-\btheta^\star) }=\Tr(Q S_k)$$ as well as $\Tr(Q)\leq d$ yields the second claimed inequality.\qed

{\it Proof of Inequality \eqref{eq.first_term_asymp}.} Set $k_*:= 2e^{\gamma + 3\kappa} d^2\log(d)/B.$ Then,
\begin{align}
    \Big(\frac{2d^2\log(d)}{Bk_*}\Big)^{\log(d)}\leq d^{-\gamma-3\kappa}
    \leq \frac 1{d^\gamma(d^2\log(d))^\kappa}\leq \frac{C'}{d^\gamma k_*^\kappa},
    \label{eq.dajadsh}
\end{align}
with $C'=(2e^{\gamma + 3\kappa}/B)^\kappa.$ If for positive $A,B,u,v,x$ with $u>v,$ we have $A/x^u\leq B/x^v,$ then also $A/y^u\leq B/y^v$ for all $y\geq x.$ Since $d>e^\kappa$ we can apply this inequality with $u=\log(d)>\kappa=v,$ proving that \eqref{eq.dajadsh} holds for all $k\geq k_*.$

\section{Proofs for the lower bounds}
\label{sec.proofs_lb}

\begin{proof}[Proof of Lemma \ref{lem.equiv_model}] To distinguish the two models, we rename the queries in \eqref{eq.zlzl'_2} into $W_1,W_1',\ldots,\linebreak W_k,W_k'$ and the query vectors into $\bw_1,\bw_1',\ldots,\bw_k,\bw_k'.$ This means that the second query model is then denoted by $\{W_1,W_1',\ldots,W_k,W_k', \bw_1,\bw_1',\ldots,\bw_k,\bw_k'\}$ with $W_\ell=Y_\ell-\bX_\ell^\top \bw_\ell$ and $W_\ell'=\bX_\ell^\top \bw_\ell'.$ 

To see the equivalence, we apply induction with respect to $k.$ The base case $k=1$ and the induction step $k\to k+1$ are similar and therefore only the latter will be discussed.

The induction step $k\to k+1$ is split in two parts. We first prove  that the first query model can be transformed into the second query model without knowledge of the parameters. To see this, choose query vectors $\bv_{k+1}=\bw_{k+1}+\bw_{k+1}'$ and $\bv_{k+1}'=\bw_{k+1}-\bw_{k+1}'.$ Query vectors can depend on previously seen queries and query vectors, such that $\bw_{k+1}, \bw_{k+1}'$ can depend on $\{W_1,W_1',\ldots,W_k,W_k', \bw_1,\bw_1',\ldots,\bw_k,\bw_k'\}$ and $\bv_{k+1},\bv_{k+1}'$ can depend on $\{Z_1,Z_1',\ldots,Z_k,Z_k',\linebreak \bv_1,\bv_1',\ldots,\bv_k,\bv_k'\}$. Since the models can be transformed into each other by the induction hypothesis, those are eligible query vectors. The corresponding queries in \eqref{eq.zlzl'} are then given by $Z_{k+1}=Y_{k+1}-\bX_{k+1}^\top (\bw_{k+1}+\bw_{k+1}')$ and $Z_{k+1}'=Y_{k+1}-\bX_{k+1}^\top (\bw_{k+1}-\bw_{k+1}').$ Now, $W_{k+1}:=\tfrac 12 (Z_{k+1}+Z_{k+1}')=Y_{k+1}-\bX_{k+1}^\top \bw_{k+1}$ and $W_{k+1}':=\tfrac 12 (Z_{k+1}'-Z_{k+1})=\bX_{k+1}^\top \bw_{k+1}'$ are then the queries from the second model. Since one can also retrieve the query vectors $\bw_{k+1}=\tfrac 12 (\bv_{k+1}+\bv_{k+1}')$ and $\bw_{k+1}'=\tfrac 12 (\bv_{k+1}-\bv_{k+1}'),$ we can transform the data from the first query model into data from the second query model $\{W_1,W_1',\ldots,W_{k+1},W_{k+1}', \bw_1,\bw_1',\ldots,\bw_{k+1},\bw_{k+1}'\},$ completing the first part of the induction step.

The induction step for the other direction is similar and omitted. 
\end{proof}

As in the previous work on lower bounds for sequential designs \cite{pmlr-v30-Shamir13, 6289365, Duchi2015}, we use a version of Assouad's lemma \cite{MR2724359}. Recall that $\theta_j$ denotes the $j$-th component of $\btheta=(\theta_1,\ldots,\theta_d)^\top.$

\begin{lem}
\label{lem.assouad}
For any estimator $\wh \btheta$ and any $\rho>0,$
\begin{align*}
    \max_{\btheta\in \btheta^\star+\{-\rho,\rho\}^d} \ \ex_{\btheta,\mV_k} \ \big[\|\widehat \btheta-\btheta\|^2\big] \geq \frac {\rho^2 d}2 
     \Big(1 - \frac{1}{d}\sum_{j=1}^d \TV\big(\mathbb{P}_{+j,\btheta^\star,\mV_k},\mathbb{P}_{-j,\btheta^\star,\mV_k}\big)\Big)
\end{align*}
with 
\begin{align*}
    \mathbb{P}_{+j,\btheta^\star,\mV_k} := \frac{1}{2^{d-1}}
    \sum_{\btheta\in \btheta^\star+\{-\rho,\rho\}^d \ \text{and} \ \theta_j=\theta_j^\star +\rho} P_{\btheta,\mV_k},
\end{align*}
and
\begin{align*}
    \mathbb{P}_{-j,\btheta^\star,\mV_k} := \frac{1}{2^{d-1}}
    \sum_{\btheta\in \btheta^\star+\{-\rho,\rho\}^d \ \text{and} \ \theta_j=\theta_j^\star -\rho} P_{\btheta,\mV_k}.
\end{align*}
\end{lem}

\begin{proof}
Let $\xi$ be either $+1$ or $-1.$ If $\btheta$ is an element of $\btheta^\star+\{-\rho,\rho\}^d$ with $\theta_j=\theta_j^\star+\xi\rho,$ then, $(\widehat \theta_j-\theta_j)^2\geq \rho^2 \mathbf{1}(\sign(\widehat \theta_j-\theta_j^\star)\neq \xi),$ where the sign function evaluated at zero is defined as $+1.$ For probability measures $P,Q$ and a measurable set $A$, we have $P(A)+Q(A^c)=1+Q(A^c)-P(A^c)\geq 1-\TV(P,Q).$ Rewriting $\|\widehat \btheta-\btheta\|^2=\sum_{j=1}^d (\widehat \theta_j-\theta_j)^2,$ and applying the lower bound $$\max_{\btheta\in \btheta^\star+\{-\rho,\rho\}^d}\ \geq \frac{1}{2^d} \ \sum_{\btheta\in \btheta^\star+\{-\rho,\rho\}^d} = \frac{1}{2^d} \ \sum_{\btheta\in \btheta^\star+\{-\rho,\rho\}^d \ \text{and} \ \theta_j=\theta_j^\star + \rho}+\frac{1}{2^d} \ \sum_{\btheta\in \btheta^\star+\{-\rho,\rho\}^d \ \text{and} \ \theta_j=\theta_j^\star - \rho},$$
we find
\begin{align}
\begin{split}
    &\max_{\btheta\in \btheta^\star+\{-\rho,\rho\}^d} \ex_{\btheta,\mV_k} \ \big[\|\widehat \btheta-\btheta\|^2\big] \\
    &= 
    \max_{\btheta\in \btheta^\star+\{-\rho,\rho\}^d}
    \sum_{j=1}^d \ex_{\btheta,\mV_k} \ \big[(\widehat \theta_j-\theta_j)^2\big] \\
    &\geq 
    \rho^2 \sum_{j=1}^d \frac 12 \mathbb{P}_{+j,\btheta^\star,\mV_k}\big(\sign(\widehat \theta_j-\theta_j^\star)\neq 1\big) + \frac 12 \mathbb{P}_{-j,\btheta^\star,\mV_k}\big(\sign(\widehat \theta_j-\theta_j^\star)\neq -1\big) \\
    &\geq \frac {\rho^2 d}2 
     \Big(1 - \frac{1}{d}\sum_{j=1}^d \TV\big(\mathbb{P}_{+j,\btheta^\star,\mV_k},\mathbb{P}_{-j,\btheta^\star,\mV_k}\big)\Big).
     \end{split}
     \label{eq.first_step_lb}
\end{align}
The right hand side does not depend anymore on the estimator $\widehat \btheta.$    
\end{proof}

\begin{proof}[Proof of Theorem \ref{thm.lb_adapt}]
For $\sigma^2=0,$ there is nothing to prove and therefore, we assume $\sigma^2>0.$ 

We work in the equivalent model \eqref{eq.zlzl'_2}. One can assume that $\bv_\ell'\neq 0$ for all $\ell=1,\ldots,k.$ Otherwise, $Z_\ell'=0$ and changing $\bv_\ell'$ to a non-zero vector yields additional information about the model. We can also assume that for any $\ell=1,\ldots,k,$ 
\begin{align}
    \|\bv_\ell'\|=1,
    \label{eq.vell'normalization}
\end{align}
as any other scaling would simply scale the observation $Z'_\ell$. Define
\begin{align}
    \tau:= \frac{R}{9\sqrt{d}} \Big(1 \wedge \frac{d\sigma}{\sqrt{k}R}\Big). 
    \label{eq.tau_def}
\end{align}
One can think of $\tau^2$ as the convergence rate of an individual component $\ex_{\btheta,\mV_k}(\wh \btheta_j-\btheta_j)^2.$

We want to apply Assouad's lower bound in Lemma \ref{lem.assouad}. For any $\btheta^\star\in B_{R/2}(0)$ and any $\btheta\in \btheta^\star+\{-\tau,\tau\}^d,$ we have $\|\btheta\|\leq \|\btheta^\star\|+R/9\leq R.$ Thus, we can
lower bound $\sup_{\btheta\in B_{R}(0)}\geq \sup_{\btheta^\star\in B_{R/2}(0)} \ \max_{\btheta\in \btheta^\star+\{-\tau,\tau\}^d}.$ Lemma \ref{lem.assouad} with $\rho$ replaced by $\tau$ gives therefore
\begin{align}
    \sup_{\btheta\in B_{R}(0)} \ \ex_{\btheta,\mV_k} \ \big[\|\widehat \btheta-\btheta\|^2\big] \geq 
    \sup_{\btheta^\star\in B_{R/2}(0)} \ \frac {\tau^2 d}2 
     \Big(1 - \frac{1}{d}\sum_{j=1}^d \TV\big(\mathbb{P}_{+j,\btheta^\star,\mV_k},\mathbb{P}_{-j,\btheta^\star,\mV_k}\big)\Big)
     \label{eq.risk_lb_step1}
\end{align}
with 
\begin{align}
    \mathbb{P}_{+j,\btheta^\star,\mV_k} := \frac{1}{2^{d-1}}
    \sum_{\btheta\in \btheta^\star+\{-\tau,\tau\}^d \ \text{and} \ \theta_j=\theta_j^\star +\tau} P_{\btheta,\mV_k},
    \label{eq.P+_def}
\end{align}
and
\begin{align}
    \mathbb{P}_{-j,\btheta^\star,\mV_k} := \frac{1}{2^{d-1}}
    \sum_{\btheta\in \btheta^\star+\{-\tau,\tau\}^d \ \text{and} \ \theta_j=\theta_j^\star -\tau} P_{\btheta,\mV_k}.
    \label{eq.P-_def}
\end{align}

Pinsker's inequality $\TV(P,Q)\leq \sqrt{\KL(P,Q)/2}$  combined with Cauchy-Schwarz inequality \linebreak $(\tfrac 1d\sum_{j=1}^d b_j)^2\leq (\sum_{j=1}^d d^{-2})(\sum_{j=1}^d b_j^2)=\tfrac 1d \sum_{j=1}^d b_j^2$ and the joint convexity of the Kullback-Leibler divergence $\KL(\lambda_1P_1+\ldots+\lambda_m P_m,\lambda_1Q_1+\ldots+\lambda_m Q_m)\leq \sum_{s=1}^m \lambda_s \KL(P_s,Q_s)$ for $\lambda_1,\ldots,\lambda_m\geq 0$ and $\sum_{s=1}^m \lambda_s=1,$ yield
\begin{align}
\begin{split}
    &\Big(\frac{1}{d}\sum_{j=1}^d \TV\big(\mathbb{P}_{+j,\btheta^\star,\mV_k},\mathbb{P}_{-j,\btheta^\star,\mV_k}\big)\Big)^2 \\
    &\leq 
    \frac{1}{d}\sum_{j=1}^d \TV\big(\mathbb{P}_{+j,\btheta^\star,\mV_k},\mathbb{P}_{-j,\btheta^\star,\mV_k}\big)^2 \\
    &\leq \frac{1}{2d} \sum_{j=1}^d 
    \KL\big(\mathbb{P}_{+j,\btheta^\star,\mV_k},\mathbb{P}_{-j,\btheta^\star,\mV_k}\big)\\
    &\leq \frac{1}{2d} \sum_{j=1}^d \frac{1}{2^{d-1}} 
    \sum_{\btheta\in \btheta^\star+\{-\tau,\tau\}^d \ \text{and} \ \theta_j=\theta_j^\star +\tau} \KL\big(P_{\btheta,\mV_k},P_{\btheta-2\tau \be_j,\mV_k}\big) \\
    &= \frac{1}{2^d d} \sum_{\btheta\in \btheta^\star+\{-\tau,\tau\}^d} \sum_{j=1}^d  
     \KL\big(P_{\btheta,\mV_k},P_{\btheta-2\tau \be_j,\mV_k}\big) \mathbf{1}\big(\theta_j=\theta_j^\star +\tau\big).
\end{split}
    \label{eq.TV_to_KL_2}
\end{align}
with $\be_j$ the $j$-th standard basis vector.

In a next step, we need to bound $\KL(P_{\btheta,\mV_k},P_{\btheta-2\tau \be_j,\mV_k}).$ The chain rule for the Kullback-Leibler divergence states that 
\begin{align*}
    \KL(P_{U,V},Q_{U,V})&=\ex_{P_V}[\KL(P_{U|V},Q_{U|V})]+\KL(P_V,Q_V)\\ &=\ex_{P_{U,V}}\big[\KL(P_{U|V},Q_{U|V})+\KL(P_V,Q_V)\big].
\end{align*}
The data are generated sequentially,
\begin{align*}
    (\bv_1,\bv_1') \to (Z_1,Z_1') \to \ldots \to (\bv_k,\bv_k') \to (Z_k,Z_k').
\end{align*}
After the query vectors $(\bv_\ell,\bv_\ell')$ are selected, the queries are given by $Z_\ell=Y_\ell-\bX_\ell^\top \bv_\ell$ and $Z_\ell'=\bX_\ell^\top \bv_\ell'.$ Therefore,
\begin{align*}
    (Z_\ell,Z_\ell')\big| 
    (Z_1,Z_1',\ldots,Z_{\ell-1},Z_{\ell-1}',\bv_1,\bv_1',\ldots,\bv_{\ell},\bv_{\ell}')
    =(Z_\ell,Z_\ell')\big|(\bv_\ell,\bv_\ell').
\end{align*}
The distribution of $(\bv_\ell,\bv_\ell')|(Z_1,Z_1',\ldots,Z_{\ell-1},Z_{\ell-1}',\bv_1,\bv_1',\ldots,\bv_{\ell-1},\bv_{\ell-1}')$ does not depend on the unknown regression vector $\btheta.$ Write $Q_\ell$ to denote this distribution. If $P_{\btheta,(Z_\ell,Z_\ell')|(\bv_\ell,\bv_\ell')}$ denotes the distribution of $(Z_\ell,Z_\ell')\big|(\bv_\ell,\bv_\ell')$ for the data generating parameter $\btheta,$ the chain rule and the arguments above yield
\begin{align}
    \begin{split}
    \KL\big(P_{\btheta,\mV_k},P_{\btheta',\mV_k}\big)
    &= \E_{\btheta,\mV_k}\Big[\sum_{\ell=1}^k \KL\big(P_{\btheta,(Z_\ell,Z_\ell')|(\bv_\ell,\bv_\ell')},P_{\btheta',(Z_\ell,Z_\ell')|(\bv_\ell,\bv_\ell')}\big)+\KL(Q_\ell,Q_\ell)\Big]\\
    &= \E_{\btheta,\mV_k}\Big[\sum_{\ell=1}^k \KL\big(P_{\btheta,(Z_\ell,Z_\ell')|(\bv_\ell,\bv_\ell')},P_{\btheta',(Z_\ell,Z_\ell')|(\bv_\ell,\bv_\ell')}\big)\Big].
    \end{split}
    \label{eq.KL_sum_k}
\end{align}
We will now apply this identity for $\btheta'=\btheta-2\tau\be_j$ to eventually derive a bound for \eqref{eq.TV_to_KL_2}.

For two centered $d$-variate normal distributions $\mN(0,\Sigma_0),$ $\mN(0,\Sigma_1),$ with $\Sigma_0,\Sigma_1>0,$ we have 
\begin{align*}
    \KL\Big(\mathcal{N}(0,\Sigma_1^2), \mathcal{N}(0,\Sigma_0^2)\Big)
    = \frac{1}{2}\bigg(\log\Big(\frac{\det(\Sigma_0)}{\det(\Sigma_1)}\Big)+
    \Tr\big(\Sigma_1\Sigma_0^{-1}\big)-d\bigg).    
\end{align*}
Consider an invertible matrix $T,$ and let $S_0:=T\Sigma_0T^\top,$ $S_1:= T\Sigma_1T^\top.$ Assume further that there are  positive definite matrices $\Lambda_0,\Lambda_1$ such that $S_0\geq \Lambda_0$ and $S_1\geq \Lambda_1.$ It is known (e.g. Theorem 6.8 in \cite{MR1691203}), that $S_0^{-1}\leq \Lambda_0^{-1},$ $S_1^{-1}\leq \Lambda_1^{-1},$ and that $\Tr(AB)=\Tr(BA)$ for square matrices $A, B.$ Using these facts, the symmetrized Kullback-Leibler divergence can be bounded as follows
\begin{align}
\begin{split}
	\KL\Big(\mathcal{N}(0,\Sigma_0^2), \mathcal{N}(0,\Sigma_1^2)\Big)
    &\leq\KL\Big(\mathcal{N}(0,\Sigma_0^2), \mathcal{N}(0,\Sigma_1^2)\Big)+\KL\Big(\mathcal{N}(0,\Sigma_1^2), \mathcal{N}(0,\Sigma_0^2)\Big) \\
    &= \frac{1}{2}\Big(
    \Tr\big(\Sigma_1\Sigma_0^{-1}\big)+\Tr\big(\Sigma_0\Sigma_1^{-1}\big)-2d\Big) \\
    &=\frac{1}{2}\Tr\big(\Sigma_1^{-1}(\Sigma_1-\Sigma_0)\Sigma_0^{-1}(\Sigma_1-\Sigma_0)\big) \\
    &= \frac{1}{2}\Tr\big((T\Sigma_1T^\top)^{-1}(T\Sigma_1T^\top-T\Sigma_0T^\top)(T\Sigma_0T^\top)^{-1}(T\Sigma_1T^\top-T\Sigma_0T^\top)\big) \\
    &= \frac{1}{2}\Tr\big(S_1^{-1}(S_1-S_0)S_0^{-1}(S_1-S_0)\big) \\
    &= \frac{1}{2}\Tr\big(S_1^{-1/2}(S_1-S_0)S_0^{-1}(S_1-S_0)S_1^{-1/2}\big) \\
    &\leq  \frac{1}{2}\Tr\big(S_1^{-1/2}(S_1-S_0)\Lambda_0^{-1}(S_1-S_0)S_1^{-1/2}\big) \\
    &= \frac{1}{2}\Tr\big(\Lambda_0^{-1/2}(S_1-S_0)S_1^{-1}(S_1-S_0)\Lambda_0^{-1/2}\big) \\
    &\leq \frac{1}{2}\Tr\big(\Lambda_0^{-1/2}(S_1-S_0)\Lambda_1^{-1}(S_1-S_0)\Lambda_0^{-1/2}\big) \\
    &= \frac{1}{2}\Tr\big(\Lambda_0^{-1}(S_1-S_0)\Lambda_1^{-1}(S_1-S_0)\big).
    \end{split}
    \label{eq.KL_formula123}
\end{align}

We now control the right hand side of \eqref{eq.KL_sum_k} for $\btheta'=\btheta-2\tau\be_j$. As we work in the equivalent model \eqref{eq.zlzl'_2} with the normalization constraint \eqref{eq.vell'normalization}, the distributions of $(Z_\ell,Z_\ell')^\top|(\bv_\ell,\bv_\ell')$ is given by \eqref{eq.zlzl'_distrib_2}. We have normalized $\bv_\ell'$ to a unit-length vector $\|\bv_\ell'\|=1.$ Thus, the distributions $P_{\btheta, (Z_\ell,Z_\ell')|(\bv_\ell,\bv_\ell')}$ and $P_{\btheta-2\tau \be_j, (Z_\ell,Z_\ell')|(\bv_\ell,\bv_\ell')}$ (as shown in \eqref{eq.zlzl'_distrib_2}) are centered normal with respective covariances
\begin{align*}
\left(
    \begin{array}{cc}
     \sigma^2+\|\btheta-\bv_\ell\|^2    &  \langle \btheta-\bv_\ell,\bv_\ell'\rangle \\
     \langle \btheta-\bv_\ell,\bv_\ell'\rangle & 1
    \end{array}
    \right)
\end{align*}
and
\begin{align*}
    \left(
    \begin{array}{cc}
     \sigma^2+\|\btheta-2\tau \be_j-\bv_\ell\|^2    &  \langle \btheta-2\tau \be_j-\bv_\ell,\bv_\ell'\rangle \\
     \langle \btheta-2\tau \be_j-\bv_\ell,\bv_\ell'\rangle & 1
    \end{array}
    \right).
\end{align*}
In a next step, we transform the model such that the subsequent analysis becomes more tractable. Choosing the transformation
\begin{align*}
    T=\left(
    \begin{array}{cc}
     1  &  -\langle \btheta-\bv_\ell,\bv_\ell'\rangle\\
    0    & 1
    \end{array}
    \right)
\end{align*}
maps $(Z_\ell,Z_\ell')^\top$ to $(Z_\ell-\langle \btheta-\bv_\ell,\bv_\ell'\rangle Z_\ell',Z_\ell')^\top$  transforming the covariances via the identities
\begin{align}
    S_0:=T \left(
    \begin{array}{cc}
     \sigma^2+\|\btheta-\bv_\ell\|^2    &  \langle \btheta-\bv_\ell,\bv_\ell'\rangle \\
     \langle \btheta-\bv_\ell,\bv_\ell'\rangle & 1
    \end{array}
    \right) T^\top 
    =  \left(
    \begin{array}{cc}
     \sigma^2 +\|\btheta-\bv_\ell -\langle \btheta-\bv_\ell, \bv_\ell'\rangle \bv_\ell'\|^2  &  0\\
     0    & 1
    \end{array}
    \right)
    \label{eq.S0_def}
\end{align}
and
\begin{align}
    \begin{split}
    S_1&:=T \left(
    \begin{array}{cc}
     \sigma^2+\|\btheta-2\tau \be_j-\bv_\ell\|^2    &  \langle \btheta-2\tau \be_j-\bv_\ell,\bv_\ell'\rangle \\
     \langle \btheta-2\tau \be_j-\bv_\ell,\bv_\ell'\rangle & 1
    \end{array}
    \right) T^\top \\
    &=  \left(
    \begin{array}{cc}
     \sigma^2 +\|\btheta-2\tau \be_j-\bv_\ell -\langle \btheta-\bv_\ell, \bv_\ell'\rangle \bv_\ell'\|^2  &  -2\tau \langle \be_j,\bv_\ell'\rangle\\
     -2\tau \langle \be_j,\bv_\ell'\rangle    & 1
    \end{array}
    \right).
    \end{split}
    \label{eq.S1_def}
\end{align}
For a symmetric $2\times 2$ matrix, the elementary inequality $2ab \geq -|2ab|\geq -2a^2-\tfrac 12 b^2$ yields for all vectors $(u_1,u_2)^\top,$
\begin{align*}
	(u_1 u_2)^\top 
	\left(
	\begin{array}{cc}
	\alpha & \beta \\
	\beta & 1
	\end{array}
	\right) 
	\left(
	\begin{array}{c}
	u_1 \\
	u_2
	\end{array}	
	\right) 
	= u_1^2 \alpha +2u_1u_2 \beta+u_2^2\geq u_1^2(\alpha - 2\beta^2) +\frac 12 u_2^2
\end{align*}
implying the matrix inequality
\begin{align}
	\left(
	\begin{array}{cc}
	\alpha & \beta \\
	\beta & 1
	\end{array}
	\right)
	\geq 
	\left(
	\begin{array}{cc}
	\alpha -2\beta^2 & 0 \\
	0 & \frac 12
	\end{array}
	\right).
    \label{eq.matrix_inequality}
\end{align}
If $k\geq d^2$, then 
\begin{align}
	\tau = \frac{R}{9\sqrt{d}}\Big(1\wedge \frac{d\sigma}{\sqrt{k} R}\Big) \leq \frac{\sigma}{9} \sqrt{\frac{d}{k}}\leq \frac{\sigma}4.
	\label{eq.tau_ub}
\end{align}
The last inequality is rather loose but sufficient to give $8\tau^2\leq \sigma^2/2.$ Together with the matrix inequality applied to $\beta=-2\tau \langle \be_j,\bv_\ell'\rangle,$ we obtain
\begin{align*}
	S_1 \geq 
	\left(
    \begin{array}{cc}
     \frac 12 \sigma^2 +\|\btheta-2\tau\be_j -\bv_\ell -\langle \btheta-\bv_\ell, \bv_\ell'\rangle \bv_\ell'\|^2  &  0\\
     0    & \frac 12
    \end{array}
    \right)
    \geq 
    \frac 12 
    \left(
    \begin{array}{cc}
     \sigma^2  &  0\\
     0    & 1
    \end{array}
    \right)
    =: \Lambda_1.
\end{align*}
Setting $\Lambda_0:=S_0,$ we can now apply \eqref{eq.KL_formula123} for these choices of $S_0,S_1,\Lambda_0, \Lambda_1.$ Observing that 
\begin{align*}
	\|\btheta-2\tau\be_j -\bv_\ell -\langle \btheta-\bv_\ell, \bv_\ell'\rangle \bv_\ell'\|^2
	- &\|\btheta -\bv_\ell -\langle \btheta-\bv_\ell, \bv_\ell'\rangle \bv_\ell'\|^2 \\
	&= 4\tau^2-4\tau \langle \be_j, \btheta -\bv_\ell -\langle \btheta-\bv_\ell, \bv_\ell'\rangle \bv_\ell'\rangle,
\end{align*}
\begin{align*}
	S_1-S_0 = 
	 \left(
    \begin{array}{cc}
      4\tau^2-4\tau \langle \be_j, \btheta -\bv_\ell -\langle \btheta-\bv_\ell, \bv_\ell'\rangle \bv_\ell'\rangle &  -2\tau \langle \be_j,\bv_\ell'\rangle\\
     -2\tau \langle \be_j,\bv_\ell'\rangle    & 0
    \end{array}
    \right),
\end{align*}
and
\begin{align}
&\frac 12 \Tr\left(\left(
    \begin{array}{cc}
     \lambda_0^{-1}  &  0\\
     0    & 1
    \end{array}
    \right)\left(
    \begin{array}{cc}
     \alpha  &  \beta\\
     \beta    & 0
    \end{array}
    \right)2\left(
    \begin{array}{cc}
     \lambda_1^{-1}  &  0\\
     0    & 1
    \end{array}
    \right)\left(
    \begin{array}{cc}
     \alpha  &  \beta\\
     \beta    & 0
    \end{array}
    \right)\right) \\
    &= 
    \Tr\left( \left(
    \begin{array}{cc}
     \alpha \lambda_0^{-1}  &  \beta\lambda_0^{-1}\\
     \beta    & 0
    \end{array}
    \right) \left(
    \begin{array}{cc}
     \alpha \lambda_1^{-1}  &  \beta\lambda_1^{-1}\\
     \beta    & 0
    \end{array}
    \right)\right) \\
    &= \Tr\left(
    \left(
    \begin{array}{cc}
     \alpha^2 \lambda_0^{-1}\lambda_1^{-1}+\beta^2 \lambda_0^{-1}  &  \alpha\beta \lambda_0^{-1}\lambda_1^{-1} \\
     \alpha\beta \lambda_1^{-1}   & \beta^2\lambda_1^{-1}
    \end{array}
    \right)
    \right) \\
    &= \frac{\alpha^2}{\lambda_0\lambda_1}  +\beta^2\Big(\frac{1}{\lambda_0}+\frac{1}{\lambda_1}\Big),
    \label{eq.trace_formula}
\end{align}
we find for the specific choices $\lambda_0=\sigma^2 +\|\btheta-\bv_\ell -\langle \btheta-\bv_\ell, \bv_\ell'\rangle \bv_\ell'\|^2,$ $\lambda_1=\sigma^2$, $\alpha^2= (4\tau^2-4\tau \langle \be_j, \btheta -\bv_\ell -\langle \btheta-\bv_\ell, \bv_\ell'\rangle \bv_\ell'\rangle)^2\leq 32 \tau^4+32\tau^2 \langle \be_j, \btheta -\bv_\ell -\langle \btheta-\bv_\ell, \bv_\ell'\rangle \bv_\ell'\rangle^2$ and $\beta= -2\tau \langle \be_j,\bv_\ell'\rangle,$ using \eqref{eq.KL_formula123}, and $\lambda_0\geq \sigma^2,$
\begin{align*}
	&\KL\big(P_{\btheta, (Z_\ell,Z_\ell')|(\bv_\ell,\bv_\ell')},P_{\btheta-2\tau \be_j, (Z_\ell,Z_\ell')|(\bv_\ell,\bv_\ell')}\big) \\
	&\leq \frac{1}{2}\Tr\big(\Lambda_0^{-1}(S_1-S_0)\Lambda_1^{-1}(S_1-S_0)\big) \\
	&= \frac{\alpha^2}{\lambda_0\lambda_1}   +\beta^2\Big(\frac{1}{\lambda_0}+\frac{1}{\lambda_1}\Big)\\
    &\leq \frac{32 \tau^4+32\tau^2 \langle \be_j, \btheta -\bv_\ell -\langle \btheta-\bv_\ell, \bv_\ell'\rangle \bv_\ell'\rangle^2}{(\sigma^2 +\|\btheta-\bv_\ell -\langle \btheta-\bv_\ell, \bv_\ell'\rangle \bv_\ell'\|^2)\sigma^2 } + 6 \frac{\tau^2}{\sigma^2} \langle \be_j,\bv_\ell'\rangle^2 \\
    &\leq 32\frac{ \tau^4}{\sigma^4}+32\frac{\tau^2 \langle \be_j, \btheta -\bv_\ell -\langle \btheta-\bv_\ell, \bv_\ell'\rangle \bv_\ell'\rangle^2 }{\sigma^2 \|\btheta-\bv_\ell -\langle \btheta-\bv_\ell, \bv_\ell'\rangle \bv_\ell'\|^2}+ 6 \frac{\tau^2}{\sigma^2} \langle \be_j,\bv_\ell'\rangle^2.
\end{align*}
Interchanging the sums over $j$ and $\ell,$ we conclude from \eqref{eq.TV_to_KL_2}, \eqref{eq.KL_sum_k}, the previous inequality, $\sum_{j=1}^d \langle \be_j,\bv_\ell'\rangle^2=\|\bv_\ell'\|^2,$ $\sum_{j=1}^d \langle \be_j, \btheta -\bv_\ell -\langle \btheta-\bv_\ell, \bv_\ell'\rangle \bv_\ell'\rangle^2=\| \btheta -\bv_\ell -\langle \btheta-\bv_\ell, \bv_\ell'\rangle \bv_\ell'\|^2,$ the normalization $\|\bv_\ell'\|=1,$ $\tau \leq \sigma \sqrt{d/(81k)}$ (as derived in \eqref{eq.tau_ub}), and $k\geq d^2,$
\begin{align*}
	\begin{split}
    &\Big(\frac{1}{d}\sum_{j=1}^d \TV(\mathbb{P}_{+j,\btheta^\star,\mV_k},\mathbb{P}_{-j,\btheta^\star,\mV_k})\Big)^2 \\
    &\leq 
    \frac{1}{d}\sum_{j=1}^d \TV(\mathbb{P}_{+j,\btheta^\star,\mV_k},\mathbb{P}_{-j,\btheta^\star,\mV_k})^2 \\
    &\leq \frac{1}{2^d d} \sum_{\btheta\in \btheta^\star+\{-\tau,\tau\}^d} \sum_{j=1}^d  
     \KL\big(P_{\btheta,\mV_k},P_{\btheta-2\tau \be_j,\mV_k}\big) \mathbf{1}\big(\theta_j=\theta_j^\star +\tau\big) \\
     &\leq \frac{1}{2^d d} \sum_{\btheta\in \btheta^\star+\{-\tau,\tau\}^d} \sum_{j=1}^d  
	\sum_{\ell=1}^k 32\frac{ \tau^4}{\sigma^4}+32\frac{\tau^2 \langle \be_j, \btheta -\bv_\ell -\langle \btheta-\bv_\ell, \bv_\ell'\rangle \bv_\ell'\rangle^2 }{\sigma^2 \|\btheta-\bv_\ell -\langle \btheta-\bv_\ell, \bv_\ell'\rangle \bv_\ell'\|^2}+ 6 \frac{\tau^2}{\sigma^2} \langle \be_j,\bv_\ell'\rangle^2 \\
	&= \frac{1}{2^d d} \sum_{\btheta\in \btheta^\star+\{-\tau,\tau\}^d} 
	\sum_{\ell=1}^k 32d\frac{ \tau^4}{\sigma^4} +38 \frac{\tau^2}{\sigma^2}  \\
	&= 32k\frac{ \tau^4}{\sigma^4} +38 \frac{\tau^2k}{\sigma^2d} \\
	&\leq \frac{32 d^2}{81^2k}+\frac{38}{81} \\
	&\leq \frac 12.
\end{split}
\end{align*}
Taking square roots, gives $1-d^{-1}\sum_{j=1}^d \TV(\mathbb{P}_{+j,\btheta^\star,\mV_k},\mathbb{P}_{-j,\btheta^\star,\mV_k})\geq 1-1/\sqrt{2}.$ Since by definition $\tau= \tfrac{R}{9\sqrt{d}} (1 \wedge \tfrac{d\sigma}{\sqrt{k}R}),$ the claimed lower bound follows from \eqref{eq.risk_lb_step1},
\begin{align*}
	 \sup_{\btheta\in \Theta} \ \ex_{\btheta,\mV_k} \ \big[\|\widehat \btheta-\btheta\|^2\big]
    &\geq \frac {\tau^2 d}2 \sup_{\btheta^\star\in B_{R/2}(0)} \ 
     \Big(1 - \frac{1}{d}\sum_{j=1}^d \TV(\mathbb{P}_{+j,\btheta^\star,\mV_k},\mathbb{P}_{-j,\btheta^\star,\mV_k})\Big) \\
     &=\frac{1}{162}\Big(1-\frac{1}{\sqrt{2}}\Big) \bigg( R^2\wedge \frac{d^2}{k}\sigma^2\bigg).
\end{align*}
\end{proof}

\begin{proof}[Proof of Theorem \ref{thm.non_adapt}]
Let
\begin{align}
    \rho=c \frac{R}{\sqrt{d}} \bigg(1 \wedge \frac{d}{\sqrt{k}}\Big(1\vee \frac{\sigma}R\Big)\bigg)
    \quad \text{for} \ \ c=2^{-8}.
    \label{eq.rho_lb_def}
\end{align}
For any $\btheta^\star\in B_{R/2}(0)$ and any $\btheta\in\btheta^\star+\{-\rho,\rho\}^d,$ we have $\|\btheta\|\leq \|\btheta^\star\|+cR\leq R.$ Thus, we can
lower bound $\sup_{\btheta\in B_R(0)}\geq \sup_{\btheta^\star\in B_{R/2}(0)} \ \max_{\btheta\in \btheta^\star+\{-\rho,\rho\}^d}.$ Lemma \ref{lem.assouad} gives
\begin{align}
    \sup_{\btheta\in B_{R}(0)} \ \ex_{\btheta,\mV_k} \ \big[\|\widehat \btheta-\btheta\|^2\big] \geq 
    \sup_{\btheta^\star\in B_{R/2}(0)} \ \frac {\rho^2 d}2 
     \Big(1 - \frac{1}{d}\sum_{j=1}^d \TV\big(\mathbb{P}_{+j,\btheta^\star,\mV_k},\mathbb{P}_{-j,\btheta^\star,\mV_k}\big)\Big)
     \label{eq.risk_lb_step1_2}
\end{align}
with $\mathbb{P}_{+j,\btheta^\star,\mV_k},\mathbb{P}_{-j,\btheta^\star,\mV_k}$ as in \eqref{eq.P+_def} and \eqref{eq.P-_def} with $\tau$ replaced by $\rho$.

Since the query vectors are deterministic, \eqref{eq.KL_sum_k} becomes now 
\begin{align}
    \KL\big(P_{\btheta,\mV_k},P_{\btheta-2\rho \be_j,\mV_k}\big)
    &= \sum_{\ell=1}^k \KL\big(P_{\btheta,(Z_\ell,Z_\ell')|(\bv_\ell,\bv_\ell')},P_{\btheta-2\rho \be_j,(Z_\ell,Z_\ell')|(\bv_\ell,\bv_\ell')}\big),
    \label{eq.KL_non_adaptive}
\end{align}
where $\btheta$ is an element in the set
$\btheta^\star +\{-\rho,\rho\}^d$ with $\theta_j=\theta_j^\star+\rho.$

Let $\bv'$ satisfy $\|\bv'\|=1$ and let $\bw, \bw'\in \{-1,0,1\}^d$. We now prove that $\|\bu\|\geq Rk^{-1/(d-1)}/4$ implies
\begin{align}
    \sigma^2+\big\|\bu + \rho \bw+ \rho\langle \bw', \bv'\rangle \bv'\big\|^2-8\rho^2
    \geq 
    \frac{\sigma^2\vee \|\bu\|^2}{8}.
    \label{eq.dsaRD}
\end{align}
To see this, use that $\|\bv'\|=1$ and $|\langle \bw', \bv'\rangle|\leq \|\bw'\|.$ Thus, triangle inequality yields $\rho \|\bw+ \langle \bw', \bv'\rangle \bv'\|\leq 2\rho \sqrt{d}.$ By \eqref{eq.rho_lb_def}, $\rho = c R d^{-1/2}(1\wedge dk^{-1/2}(1\vee (\sigma/R)).$ 

If $k\leq 1/(16c)^{d-1},$ we have $\|\bu\|\geq 4Rc\geq 4\rho \sqrt{d}.$ Moreover, if $k>1/(16c)^{d-1}$ and $\sigma\leq R,$ we can write $k=(a/16c)^{d-1}$ with $a>1.$ Since $y\leq e^{y-1}$ and $c\leq 1/(16e^2),$ we have $d\leq (e^2)^{(d-1)/2}\leq 1/(16c)^{(d-1)/2}.$ Recalling that $d\geq 3,$ we find again $\|\bu\|\geq Rk^{-1/(d-1)}/4=4cR/a\geq 4cRd/(a^{(d-1)/2}(1/16c)^{(d-1)/2})= 4cRd/\sqrt{k}\geq 4\rho \sqrt{d}.$ Thus in both of the previous cases $k\leq 1/(16c)^{d-1}$ and $\{k>1/(16c)^{d-1}\}  \cap \{\sigma\leq R\},$ $\|\bu\|\geq 4 \rho\sqrt{d}\geq 2\rho \|\bw+ \langle \bw', \bv'\rangle \bv'\|$ and therefore by triangle inequality, $ \|\bu+\rho\bw+ \rho\langle \bw', \bv'\rangle \bv'\|\geq \|\bu\|/2.$ Since for $d\geq 4,$ also $\|\bu\|^2/8\geq 2\rho^2 d\geq 8\rho^2,$ \eqref{eq.dsaRD} follows for these cases.

The remaining case $k>1/(16c)^{d-1}$ and $\sigma>R,$ yields $\sigma/4\geq 2c dk^{-1/2}\sigma \geq 2\rho\sqrt{d}.$ Together with the elementary inequality $a^2+b^2\geq \tfrac 12 (a+b)^2$, triangle inequality and $\rho \|\bw+ \langle \bw', \bv'\rangle \bv'\|\leq 2\rho \sqrt{d}$, we get $\sigma+\big\|\bu + \rho \bw+ \rho\langle \bw', \bv'\rangle \bv'\big\|\geq \sigma+\big\|\bu\big\|-2\rho \sqrt{d}\geq (3/4)\sigma+\big\|\bu\big\|$ and  
\begin{align*}
    \sigma^2+\big\|\bu + \rho \bw + \rho\langle \bw', \bv'\rangle \bv'\big\|^2
    &\geq \frac{1}{2}\Big(\sigma+\big\|\bu + \rho \bw+\rho\langle \bw', \bv'\rangle \bv'\big\|\Big)^2 \geq \frac{\sigma^2\vee \|\bu\|^2}{4}.
\end{align*}
Moreover $\sigma^2/8\geq 8\rho^2 d\geq 8\rho^2.$ Thus,  \eqref{eq.dsaRD} holds in all cases.

Let
\begin{align*}
    S_0=
    \left(
    \begin{array}{cc}
     \sigma^2 +\|\btheta-\bv_\ell -\langle \btheta-\bv_\ell, \bv_\ell'\rangle \bv_\ell'\|^2  &  0\\
     0    & 1
    \end{array}
    \right)
\end{align*}
be as in \eqref{eq.S0_def}.

For the next arguments, we will assume that 
\begin{align}
    \big\|\btheta^\star-\bv_\ell+\langle \btheta^\star-\bv_\ell,\bv_\ell'\rangle\big\|
    \geq \frac R4 k^{-1/(d-1)}.
    \label{eq.sep_condition_lb}
\end{align}
Note that $\btheta\in \btheta^\star+\{-\rho,\rho\}^d.$ Thus, applying \eqref{eq.dsaRD} with $\bu=\btheta^\star-\bv_\ell+\langle \btheta^\star-\bv_\ell,\bv_\ell'\rangle$, $\rho\bw=\btheta-\btheta^\star$, we have that $S_0\geq \wt \Lambda_0$ with 
\begin{align*}
    \wt \Lambda_0=
    \left(
    \begin{array}{cc}
     \frac{\sigma^2 \vee \|\btheta^\star-\bv_\ell -\langle \btheta^\star-\bv_\ell, \bv_\ell'\rangle \bv_\ell'\|^2}8  &  0\\
     0    & 1
    \end{array}
    \right).
\end{align*}
Similarly, let
\begin{align*}
    S_1=  \left(
    \begin{array}{cc}
     \sigma^2 +\|\btheta-2\rho \be_j-\bv_\ell -\langle \btheta-\bv_\ell, \bv_\ell'\rangle \bv_\ell'\|^2  &  -2\rho \langle \be_j,\bv_\ell'\rangle\\
     -2\rho \langle \be_j,\bv_\ell'\rangle    & 1
    \end{array}
    \right)
\end{align*}
be as in \eqref{eq.S1_def}. Using the matrix inequality \eqref{eq.matrix_inequality}, $\langle \be_j,\bv_\ell'\rangle\leq \|\be_j\| \|\bv_\ell'\|=1$, and applying \eqref{eq.dsaRD} with $\bw=\btheta-\btheta^\star-2\rho \be_j$ and $\bw'=\btheta-\btheta^\star,$ we have
\begin{align*}
    S_1
    &\geq   \left(
    \begin{array}{cc}
     \sigma^2 +\|\btheta-2\rho \be_j-\bv_\ell -\langle \btheta-\bv_\ell, \bv_\ell'\rangle \bv_\ell'\|^2 -8\rho^2  &  0\\
     0   & \tfrac 12
    \end{array}
    \right) \\
    &\geq 
     \left(
    \begin{array}{cc}
     \frac{\sigma^2 \vee \|\btheta^\star-\bv_\ell -\langle \btheta^\star-\bv_\ell, \bv_\ell'\rangle \bv_\ell'\|^2}8  &  0\\
     0   & \tfrac 12
    \end{array}
    \right).
\end{align*}

Recall that $\btheta$ is an element in the set
$\btheta^\star +\{-\rho,\rho\}^d$ with $\theta_j=\theta_j^\star+\rho.$ This means that $\langle \btheta - \rho \be_j- \btheta^\star, \be_j\rangle =0$ and thus
\begin{align*}
	\|\btheta-2\rho\be_j -\bv_\ell -\langle \btheta-\bv_\ell, \bv_\ell'\rangle \bv_\ell'\|^2
	- &\|\btheta -\bv_\ell -\langle \btheta-\bv_\ell, \bv_\ell'\rangle \bv_\ell'\|^2 \\
	&= -4\rho \langle \be_j, \btheta -\rho \be_j-\bv_\ell -\langle \btheta-\bv_\ell, \bv_\ell'\rangle \bv_\ell'\rangle \\
        &= -4\rho \langle \be_j, \btheta^\star -\langle \btheta-\bv_\ell, \bv_\ell'\rangle \bv_\ell'\rangle.
\end{align*}
This implies
\begin{align*}
	S_1-S_0 = 
	 \left(
    \begin{array}{cc}
      -4\rho \langle \be_j, \btheta^\star -\langle \btheta-\bv_\ell, \bv_\ell'\rangle \bv_\ell'\rangle &  -2\rho \langle \be_j,\bv_\ell'\rangle\\
     -2\rho \langle \be_j,\bv_\ell'\rangle    & 0
    \end{array}
    \right)
    =
    \left(
    \begin{array}{cc}
      \alpha &  \beta\\
     \beta & 0
    \end{array}
    \right).
\end{align*}
Combining \eqref{eq.KL_formula123} and \eqref{eq.trace_formula} with $\lambda_0 = \lambda_1/2 = \tfrac 18(\sigma^2 \vee \|\btheta^\star-\bv_\ell -\langle \btheta^\star-\bv_\ell, \bv_\ell'\rangle \bv_\ell'\|^2),$ and $\alpha,\beta$ as above, we obtain 
\begin{align*}
    \KL\big(P_{\btheta,(Z_\ell,Z_\ell')|(\bv_\ell,\bv_\ell')},P_{\btheta-2\rho\be_j,(Z_\ell,Z_\ell')|(\bv_\ell,\bv_\ell')}\big)
    \leq 
    &32\frac{(-4\rho \langle \be_j, \btheta^\star -\langle \btheta-\bv_\ell, \bv_\ell'\rangle \bv_\ell'\rangle)^2}{ (\sigma^2 \vee \|\btheta^\star-\bv_\ell -\langle \btheta^\star-\bv_\ell, \bv_\ell'\rangle \bv_\ell'\|^2)^2} \\
    &+ 12\frac{(2\rho \langle \be_j,\bv_\ell'\rangle)^2}{\sigma^2 \vee \|\btheta^\star-\bv_\ell -\langle \btheta^\star-\bv_\ell, \bv_\ell'\rangle \bv_\ell'\|^2}.
\end{align*}
Since $\langle \btheta-\btheta^\star,\bv_\ell'\rangle^2\leq \|\btheta-\btheta^\star\|^2 \|\bv_\ell'\|^2 =\rho^2 d,$ we find using $(a+b)^2\leq 2a^2+2b^2,$
\begin{align*}
    \langle \be_j, \btheta^\star -\langle \btheta-\bv_\ell, \bv_\ell'\rangle \bv_\ell'\rangle^2
    \leq 2 \langle \be_j, \btheta^\star -\langle \btheta^\star-\bv_\ell, \bv_\ell'\rangle \bv_\ell'\rangle^2
    +2\rho^2 d \langle \be_j,\bv_\ell'\rangle^2.
\end{align*}
Since $\be_1,\ldots,\be_d$ form an orthonormal basis, $\sum_{j=1}^d \langle \be_j,\ba\rangle^2 =\|\ba\|^2$ for all $d$-dimensional vectors $\ba.$ Thus taking the sum over $j$ yields
\begin{align*}
    &\sum_{j=1}^d \KL\big(P_{\btheta,(Z_\ell,Z_\ell')|(\bv_\ell,\bv_\ell')},P_{\btheta-2\rho\be_j,(Z_\ell,Z_\ell')|(\bv_\ell,\bv_\ell')}\big) \\
    &\leq 
    (2^8+48)\frac{\rho^2}
    {\sigma^2 \vee \|\btheta^\star-\bv_\ell -\langle \btheta^\star-\bv_\ell, \bv_\ell'\rangle \bv_\ell'\|^2}
    + 2^8 \frac{\rho^4d}{ (\sigma^2 \vee \|\btheta^\star-\bv_\ell -\langle \btheta^\star-\bv_\ell, \bv_\ell'\rangle \bv_\ell'\|^2)^2}.
\end{align*}
Recall that by \eqref{eq.KL_non_adaptive}, $\KL\big(P_{\btheta,\mV_k},P_{\btheta-2\rho \be_j,\mV_k}\big)= \sum_{\ell=1}^k \KL\big(P_{\btheta,(Z_\ell,Z_\ell')|(\bv_\ell,\bv_\ell')},P_{\btheta-2\rho \be_j,(Z_\ell,Z_\ell')|(\bv_\ell,\bv_\ell')}\big).$ Hence, 
\begin{align*}
    &\sum_{j=1}^d \KL\big(P_{\btheta,\mV_k},P_{\btheta-2\rho \be_j,\mV_k}\big) \\
    &\leq 
     \sum_{\ell=1}^k \frac{176\rho^2}
    {\sigma^2 \vee \|\btheta^\star-\bv_\ell -\langle \btheta^\star-\bv_\ell, \bv_\ell'\rangle \bv_\ell'\|^2}
    + \sum_{\ell=1}^k  \frac{2^8\rho^4d}{ (\sigma^2 \vee \|\btheta^\star-\bv_\ell -\langle \btheta^\star-\bv_\ell, \bv_\ell'\rangle \bv_\ell'\|^2)^2}.
\end{align*}
The right hand side does not depend on $\btheta$ anymore. Combined with \eqref{eq.TV_to_KL_2},
\begin{align}
    \begin{split}
    &\Big(\frac{1}{d}\sum_{j=1}^d \TV\big(\mathbb{P}_{+j,\btheta^\star,\mV_k},\mathbb{P}_{-j,\btheta^\star,\mV_k}\big)\Big)^2 \\
    &\leq \frac{1}{2^d d} \sum_{\btheta\in \btheta^\star+\{-\tau,\tau\}^d} \sum_{j=1}^d  
     \KL\big(P_{\btheta,\mV_k},P_{\btheta-2\tau \be_j,\mV_k}\big) \mathbf{1}\big(\theta_j=\theta_j^\star +\tau\big)\\
     &\leq \sum_{\ell=1}^k \frac{88\rho^2}
    {d(\sigma^2 \vee \|\btheta^\star-\bv_\ell -\langle \btheta^\star-\bv_\ell, \bv_\ell'\rangle \bv_\ell'\|^2)}
    + \sum_{\ell=1}^k  \frac{2^7\rho^4}{ (\sigma^2 \vee \|\btheta^\star-\bv_\ell -\langle \btheta^\star-\bv_\ell, \bv_\ell'\rangle \bv_\ell'\|^2)^2}.
    \end{split}
    \label{eq.78t77}
\end{align}
Notice that this bound assumes the inequality \eqref{eq.sep_condition_lb}.

If $\|\bv'\|=1,$ then $\|\bu-\bv-\langle \bu-\bv,\bv'\rangle \bv'\|=\inf_{\gamma \in \mathbb{R}}\|\bu-\bv-\gamma\bv'\|.$ Thus, the set $T_{\bv,\bv'}(r):=\{\bu:\|\bu-\bv-\langle \bu-\bv,\bv'\rangle \bv'\|\leq r\}$ is a tube with radius $r>0.$ The Euclidean ball with radius $r'>0$ and center $\bw$ is denoted by $B_{r'}(\bw):=\{\by:\|\by-\bw\|\leq r'\}.$ If $\bu\in B_{r'}(0),$ then $\langle -\bv,\bv'\rangle -r'\leq \langle \bu-\bv,\bv'\rangle \leq \langle -\bv,\bv'\rangle +r'.$ Thus, 
\begin{align*}
    H_{\bv,\bv'}(r,r'):= T_{\bv,\bv'}(r) \cap B_{r'}(0)
    \subseteq 
    \Big\{\bu:\inf_{\langle -\bv,\bv'\rangle -r'\leq \gamma \leq \langle -\bv,\bv'\rangle +r'}\|\bu-\bv-\gamma\bv'\|\leq r\Big\}.
\end{align*}
Let $\Gamma(\cdot)$ be the Gamma function. The volume formula for a tube (e.g. (2.3) in \cite{Johnstone89}) says that the volume of the right hand side is bounded by $2r'c_{d-1}r^{d-1}+c_dr^d$ with $c_p:=\pi^{p/2}/\Gamma(p/2+1)$ the volume of a $p$-dimensional unit ball. By (41) in \cite{MR3339476}, $\Gamma(d/2+1)/\Gamma(d/2+1/2)\leq \sqrt{d/2+1/2}.$ Therefore $c_{d-1}\leq \sqrt{\pi d} c_d.$ Moreover, $d\geq 5$ implies $1+2\sqrt{\pi} \sqrt{d/2+1/2} \leq 2^{d-2}.$ Thus, $2c_{d-1}+c_d\leq (1+2\sqrt{\pi} \sqrt{d/2+1/2})c_d\leq 2^{d-2} c_d.$ Together with the previous display, this means that
\begin{align*}
    \Vol\big(H_{\bv,\bv'}(r,r')\big)
    \leq 2^{d-2} r'c_dr^{d-1}, \quad \text{whenever} \ \ r\leq r'.
\end{align*}
with $\Vol$ the volume (Lebesgue measure) in $d$-dimensions. Using again that $d\geq 4,$ this means that for any $\bv_1,\bv_1',\ldots,\bv_k,\bv_k',$ and any $b>0,$
\begin{align*}
    \Vol\bigg(H_{\bv_\ell,\bv_\ell'}\Big(\frac{bR}4 k^{-1/(d-1)},\frac R 2\Big)\bigg)
    &\leq  2^{d-3} R c_d b^{d-1} R^{d-1} k^{-1} 4^{1-d} \\
    &= \frac{b^{d-1}}{2k}
    \Big(\frac R 2\Big)^d c_d \\
    &= \frac{b^{d-1}}{2k} \Vol\big(B_{R/2}(0)\big)
\end{align*}
and
\begin{align*}
    \Vol\bigg(\bigcup_{\ell=1}^k H_{\bv_\ell,\bv_\ell'}\Big(\frac R 4 k^{-1/(d-1)},\frac R 2\Big)\bigg)
    &\leq 
    \sum_{\ell=1}^k \Vol\bigg(H_{\bv_\ell,\bv_\ell'}\Big(\frac R4 k^{-1/(d-1)},\frac R 2\Big)\bigg) \\
    &\leq \frac12 \Vol\big(B_{R/2}(0)\big).
\end{align*}
The latter implies that 
\begin{align*}
    \Vol\bigg(B_{R/2}(0) \setminus \bigcup_{\ell=1}^k H_{\bv_\ell,\bv_\ell'}\Big(\frac R8 k^{-1/(d-1)},\frac R 2\Big)\bigg)
    &\geq \Vol\big(B_{R/2}(0)\big)-\frac12 \Vol\big(B_{R/2}(0)\big) \\
    &=
    \frac12 \Vol\big(B_{R/2}(0)\big).
\end{align*}
Now we define a probability measure on the ball $B_{R/2}(0)$ by 
\begin{align*}
    \nu(A) = \frac{\Vol\big(A \cap \big(B_{R/2}(0) \setminus \cup_{\ell=1}^k H_{\bv_\ell,\bv_\ell'}\big(Rk^{-1/(d-1)}/4,R/2\big)\big)\big)}{\Vol\big(B_{R/2}(0) \setminus \cup_{\ell=1}^k H_{\bv_\ell,\bv_\ell'}\big(Rk^{-1/(d-1)}/4,R/2\big)\big)}.
\end{align*}
The distribution ensures that the inequality \eqref{eq.sep_condition_lb} holds for all $\ell=1,\ldots,k$ with probability one. Thus,
we can apply \eqref{eq.78t77} and therefore lower bound $\inf_{\btheta^\star \in B_{R/2}(0)}$ by an average with respect to the probability measure $\nu,$
\begin{align}
    &\inf_{\btheta^\star \in B_{R/2}(0)}\Big(\frac{1}{d}\sum_{j=1}^d \TV\big(\mathbb{P}_{+j,\btheta^\star,\mV_k},\mathbb{P}_{-j,\btheta^\star,\mV_k}\big)\Big)^2     \notag  \\
     &\leq \int \bigg(\sum_{\ell=1}^k \frac{88\rho^2}
    {d(\sigma^2 \vee \|\btheta-\bv_\ell -\langle \btheta-\bv_\ell, \bv_\ell'\rangle \bv_\ell'\|^2)}
    + \sum_{\ell=1}^k  \frac{2^7\rho^4}{ (\sigma^2 \vee \|\btheta-\bv_\ell -\langle \btheta-\bv_\ell, \bv_\ell'\rangle \bv_\ell'\|^2)^2}\bigg) \, d\nu(\btheta)\notag \\
    &\leq \sum_{\ell=1}^k \frac{88\rho^2}{d}
    \bigg(\frac 1{\sigma^2} \wedge \int \frac 1
    {\|\btheta-\bv_\ell -\langle \btheta-\bv_\ell, \bv_\ell'\rangle \bv_\ell'\|^2}\bigg) \, d\nu(\btheta) \notag  \\
    &\quad\quad
    +  2^7\rho^4\bigg( \frac{1}{\sigma^4} \wedge  \int \frac{1}{ \|\btheta-\bv_\ell -\langle \btheta-\bv_\ell, \bv_\ell'\rangle \bv_\ell'\|^4}\bigg) \, d\nu(\btheta).\label{eq.2345}
\end{align}
The sum $\sum_{0\leq t\leq \log_2(k^{1/(d-1)})-1}$ is set to zero if $k<2^{d-1}.$ Using that for $d\geq 6,$ one has $2^{d-5}-1\geq \tfrac 12 2^{d-5},$ we obtain $\sum_{t=0}^q 2^{t(d-5)}=(2^{(q+1)(d-5)}-1)/(2^{d-5}-1)\leq 2 \cdot 2^{q(d-5)}.$ Applying this with $q=\lfloor\log_2(k^{1/(d-1)})-1\rfloor$ and using that $k^{-4/(d-1)}k=k^{(d-5)/(d-1)}$ yields
\begin{align}
    &\int \frac{1}{\|\btheta -\bv_\ell-\langle \btheta -\bv_\ell,\bv_\ell'\rangle \bv_\ell'\|^4} d\nu(\btheta) \notag \\
    &\leq \int_{\|\btheta -\bv_\ell-\langle \btheta -\bv_\ell,\bv_\ell'\rangle \bv_\ell'\|^2\geq R/8}^\infty  \frac{1}{\|\btheta -\bv_\ell-\langle \btheta -\bv_\ell,\bv_\ell'\rangle \bv_\ell'\|^4} \, d\nu(\btheta) \notag \\
    &\quad + \frac{2}{\Vol\big(B_{R/2}(0)\big)}
    \sum_{0\leq t\leq \log_2(k^{1/(d-1)})-1} \notag \\
    &\quad \quad 
    \int_{2^t R k^{-1/(d-1)}/4\leq \|\btheta -\bv_\ell-\langle \btheta -\bv_\ell,\bv_\ell'\rangle \bv_\ell'\|^2\leq 2^{t+1} R k^{-1/(d-1)}/4}
    \frac{1}{\|\btheta -\bv_\ell-\langle \btheta -\bv_\ell,\bv_\ell'\rangle \bv_\ell'\|^4} \, d\nu(\btheta) \notag \\
    &\leq 
    \frac{64}{R^4}
    + \frac{2}{\Vol\big(B_{R/2}(0)\big)}
    \sum_{0\leq t\leq \log_2(k^{1/(d-1)})-1}
    \frac{\Vol\big( H_{\bv_\ell,\bv_\ell'}\big(2^{t+1}Rk^{-1/(d-1)}/4,R/2\big)\big)}{(2^{2t} R^2 k^{-2/(d-1)}/16)^2} \notag  \\
    &=\frac{64}{R^4}
    + \frac{2}{\Vol\big(B_{R/2}(0)\big)}
    \sum_{0\leq t\leq \log_2(k^{1/(d-1)})-1}
    \frac{2^8}{2^{4t} R^4 k^{-4/(d-1)}} \frac{2^{(t+1)(d-1)}}{2k} \Vol\big(B_{R/2}(0)\big) \notag  \\
    &\leq 
    \frac{64}{R^4}
    + \frac{2^8\cdot 2^5}{R^4} \notag  \\
    &\leq \frac{2^{14}}{R^4}.\label{eq.dbjasdvj66}
\end{align}
Since $\nu$ is a probability measure, Jensen's inequality gives
\begin{align*}
    \int \frac{1}{\|\btheta -\bv_\ell-\langle \btheta -\bv_\ell,\bv_\ell'\rangle \bv_\ell'\|^2} d\nu(\btheta)
    \leq 
    \Big(\int \frac{1}{\|\btheta -\bv_\ell-\langle \btheta -\bv_\ell,\bv_\ell'\rangle \bv_\ell'\|^4} d\nu(\btheta)\Big)^{1/2}
    \leq 
    \frac{2^{7}}{R^2}.
\end{align*}
Recall that $\rho= cRd^{-1/2}(1\wedge dk^{-1/2}(1\vee \sigma/R)).$ Thus, $\rho^2\leq c^2 d k^{-1}(R\vee \sigma)^2,$ $\rho^2\leq c^2R^2d^{-1}$, and $\rho^4\leq \rho^2\rho^2\leq c^4 k^{-1}(R\vee \sigma)^2 R^2.$ Combined with \eqref{eq.2345}, \eqref{eq.dbjasdvj66}, and $c=2^{-8},$
\begin{align*}
     \begin{split}
    &\inf_{\btheta^\star \in B_{R/2}(0)}\Big(\frac{1}{d}\sum_{j=1}^d \TV\big(\mathbb{P}_{+j,\btheta^\star,\mV_k},\mathbb{P}_{-j,\btheta^\star,\mV_k}\big)\Big)^2 \leq \frac{88\cdot 2^7 k\rho^2}{d(\sigma^2\vee R^2)}+ \frac{2^{21}k\rho^4}{\sigma^4 \vee R^4}
    \leq 88\cdot 2^7 c^2+ 2^{21}c^4
    \leq \frac 14.
    \end{split}    
\end{align*}
Using \eqref{eq.risk_lb_step1_2}, we find that 
\begin{align*}
    \sup_{\btheta\in B_{R}(0)} \ \ex_{\btheta,\mV_k} \ \big[\|\widehat \btheta-\btheta\|^2\big] \geq 
    \frac {\rho^2 d}4
    = 2^{-18} R^2\Big(1\wedge \frac{d^2}{k}(R\vee \sigma)^2\Big).
\end{align*}
\end{proof}

\begin{proof}[Proof of Corollary \ref{cor.gap}]
Using $\btheta_0=0$ and Lemma \ref{lem.spect_norm_of_vec_prod} with $\ba=\bb=\btheta^\star,$ we obtain that $\|S_0\|=\|\btheta^\star(\btheta^\star)^\top\|=\|\btheta^\star\|^2\leq d.$ Thus Theorem \ref{thm.rate_ub} combined with \eqref{eq.first_term_asymp} for $\kappa=1$ and $\gamma=2$ yields the rate $d^2\log^2(d)/k$ for the upper bound in the adaptive query setting. 

The minimax rate in the non-adaptive query setting is $R^2\wedge \tfrac{d^2}k (R\vee \sigma)^2.$ Using that $R=\sqrt{d}\geq \sigma,$ this becomes $d^2\wedge d^3/k.$ Since $k\gtrsim d^2\log(d)$ the second term always dominates and the minimax rate is $d^3/k.$
\end{proof}

\section*{Acknowledgement}

We are grateful to Chao Gao for fruitful discussions. The research has been supported by the NWO Vidi grant
VI.Vidi.192.021.

\bibliographystyle{acm}       
\bibliography{bib}           
\end{document}